\documentclass[11pt, a4paper]{article}

\usepackage{amsmath,amssymb,amsthm,bm}
\usepackage{xcolor,booktabs}
\usepackage[colorlinks=true,anchorcolor=blue,filecolor=blue,linkcolor=red,urlcolor=blue,citecolor=blue]{hyperref}

\usepackage{tikz}
\usetikzlibrary{arrows.meta, calc, positioning, shapes.geometric}
\usepackage{geometry}
\usepackage{float} %%Fix the position of Figures and Tables

\AtBeginDocument{
  \setlength{\abovedisplayskip}{6pt plus 2pt minus 3pt}
  \setlength{\belowdisplayskip}{6pt plus 2pt minus 3pt}
  \setlength{\abovedisplayshortskip}{3pt plus 2pt minus 1pt}
  \setlength{\belowdisplayshortskip}{5pt plus 2pt minus 2pt}
}

\allowdisplaybreaks
\DeclareMathOperator{\bk}{bk}
\DeclareMathOperator{\Sat}{Sat}
\newcommand{\tF}{\widetilde f}
\newcommand{\by}{\bm{y}}

\newtheorem{theorem}{Theorem}[section]
\newtheorem{lemma}[theorem]{Lemma}
\newtheorem{proposition}[theorem]{Proposition}
\newtheorem{corollary}[theorem]{Corollary}

\renewcommand{\le}{\leqslant}

\renewcommand{\ge}{\geqslant}

\title{Counting $K_4$-saturating non-edges in Nosal graphs}

\author{Yongtao Li\thanks{Yau Mathematical Sciences Center, Tsinghua University, Beijing, China. 
E-mail: \texttt{ytli0921@hnu.edu.cn} } 
 }

\date{}

\hypersetup{pdftitle={Counting K4-saturating non-edges in Nosal graphs (v43)},pdfauthor={Yongtao Li}}

\begin{document}

\maketitle

\vspace{-0.7cm}

\begin{abstract}
A non-edge of a $K_4$-free graph is $K_4$-saturating if adding it
creates a copy of $K_4$, and $f(G)$ counts these non-edges. Erd\H{o}s
and Tuza showed that the minimum of $f$ over $n$-vertex $K_4$-free
graphs jumps from zero to a quadratic quantity once the number of
edges exceeds $\lfloor n^2/4\rfloor$. Balogh and Liu [{J. Combin. Theory Ser. B} {109} (2014)] determined
the asymptotic size of this jump. In this paper, we study the jump phenomenon at Nosal's spectral threshold $\lambda(G)>\sqrt m$, where $\lambda (G)$ is the spectral radius and $m$ is the number of edges.
This spectral threshold does not involve the order, and its boundary contains
unbalanced graphs far below the Mantel threshold. 
 We prove that every $K_4$-free graph $G$ with
$\lambda(G)>\sqrt m$ satisfies $f(G)+ \frac{1}{2}q(G)> \frac{1}{9} \lambda(G)^2$, where
$q(G)$ counts the vertices in triangles. The constant $\frac{1}{9}$ is
best possible, and the triangle-free case of this inequality reduces to Nosal's theorem. 
Consequently, the minimum
of $f$ over $K_4$-free graphs with $m$ edges jumps from
$0$ to $\frac19m+O(\sqrt m\,)$ at the Nosal threshold, and
the error order $\sqrt m$ is best possible.  
For large
$\lambda(G)$, we also prove
$f(G)>\frac19\lambda(G)^2-\frac13\lambda(G)$, and graphs on
the boundary $\lambda(G)^2=m$ attain this value. In the
extremal graphs, the three edges of every triangle have
booksizes in the ratio $4:1:1$, so these graphs are 
different from those of Balogh and Liu. We also show that
the spectral Mantel condition $\lambda(G)>\lambda(T_{n,2})$
forces asymptotically fewer $K_4$-saturating non-edges
than the condition $e(G)>\lfloor n^2/4\rfloor$. Our proof
transfers the Balogh--Liu theorem through a weighted
blow-up and replaces the choice of an extremal edge by a
zero-sum comparison over all edges of a neighborhood. We
expect this comparison to apply to other counts that
involve several pages of a triangle not only the largest
one. At the end of the paper, we propose several open
problems for interested readers.
\end{abstract}



%%\tableofcontents

\section{Introduction}\label{sec-1}

Mantel~\cite{Man1907} proved that every triangle-free graph on $n$
vertices has at most $\lfloor n^2/4\rfloor$ edges, and supersaturation
results describe what happens beyond this threshold. 
As a starting point, 
Rademacher
(see, e.g., Erd\H{o}s~\cite{Erd1955,erdos1962,Erdos1964})
proved that $\lfloor n^2/4\rfloor+1$ edges force at least
$\lfloor n/2\rfloor$ triangles,  
so the number of triangles jumps from zero to a linear quantity. 
Erd\H{o}s~\cite{erdos1962} extended
this bound to $k\lfloor n/2\rfloor$ triangles for graphs with at
least $\lfloor n^2/4\rfloor+k$ edges and $1\le k\le cn$, where $c>0$ is a small
constant. Lov\'asz and Simonovits~\cite{LS1975} proved the same
bound for every $1\le k<n/2$.  
Erd\H{o}s~\cite{Erd1962b} and Lov\'asz and Simonovits~\cite{LS1983}
proved analogous supersaturation for cliques. 
Mubayi~\cite{Mub2010}, and Pikhurko and
Yilma~\cite{PY2017} extended these results to color-critical graphs, 
which are graphs with an
edge whose deletion lowers the chromatic number. 
Ma and Yuan~\cite{MY2025} extended these results beyond color-critical graphs, 
for example to Kneser graphs and disjoint unions of cliques. 
For the whole range of edge densities, 
Razborov~\cite{razborov} and Nikiforov \cite{Nik2011Trans} 
determined the asymptotic minimum number of
triangles and of copies of $K_4$, respectively. 
In a breakthrough, 
Reiher~\cite{reiher} proved the clique density theorem, 
which solved a conjecture of Lov\'asz and Simonovits~\cite{LS1983}. 
Liu, Pikhurko and Staden~\cite{LPS2020} found the
exact minimum number of triangles for large graphs whose edge density
is bounded away from $1$. 
 All these results count copies of a fixed graph.
In this paper, we study a jump of another kind: we count
non-edges whose addition creates a copy of $K_4$, rather
than copies of a subgraph. 

\smallskip 
Let $G$ be a $K_4$-free graph on $n$ vertices with $m=e(G)$ edges. A
non-edge $uv$ of $G$ is \emph{$K_4$-saturating} if $G+uv$ contains a
copy of $K_4$, that is, if the common neighborhood $N(u)\cap N(v)$
contains an edge. We write $\Sat(G)$ for the set of such non-edges
and put $f(G)=|\Sat(G)|$. The literature calls these pairs
$K_4$-saturating edges~\cite{BL2014,HMMY2023}, but we say \emph{non-edges}
because they are absent from $G$.
The balanced complete bipartite graph $T_{n,2}$ has
$\lfloor n^2/4\rfloor$ edges and satisfies $f(T_{n,2})=0$. Let $f(n,m)$
denote the minimum of $f(G)$ over all $n$-vertex $K_4$-free graphs $G$
with $m$ edges. Erd\H{o}s and Tuza (see~\cite{erdos1990}) proved that
one more edge than $\lfloor n^2/4\rfloor$ already forces $\Omega(n^2)$ $K_4$-saturating non-edges. 
A natural construction is the following: add a vertex $v$ to the complete bipartite graph
 $K_{\lceil n/2\rceil -1 ,\lfloor n/2\rfloor}$, and make $v$ adjacent to roughly half of the vertices in each partite set. 
This graph has about $\frac{1}{16}n^2$ $K_4$-saturating non-edges,
and Erd\H{o}s and Tuza conjectured that
$f(n,\lfloor n^2/4\rfloor+1)=(\frac1{16}+o(1))n^2$. 
In 2014, Balogh and Liu~\cite{BL2014} disproved the Erd\H{o}s--Tuza conjecture. They proved that
$f(n,\lfloor n^2/4\rfloor+1)\ge\frac{2}{33}n^2-\frac{3}{11}n$ for
every $n\ge73$, and they constructed graphs with exactly
$\frac{2}{33}n^2-\frac{7}{33}n$ saturating non-edges whenever $66$
divides $n$. In particular,  
\begin{equation}\label{e-1}
f(n, \lfloor n^2/4\rfloor +1) = \frac{2}{33}n^2 + O(n) .
\end{equation}
Thus, the count jumps from zero to about $\frac{2}{33}n^2$ as soon as $e(G)> \lfloor n^2/4\rfloor$. 
 Let $\mathrm{ex}(n,K_{r+1})$ denote the Tur\'an number, that is, the
largest number of edges of a $K_{r+1}$-free graph on $n$ vertices. He,
Ma, Ma and Ye~\cite{HMMY2023} extended~\eqref{e-1} to
$K_{r+2}$-saturating non-edges in $K_{r+2}$-free graphs with
$\mathrm{ex}(n,K_{r+1})+1$ edges, which confirmed a conjecture of Balogh
and Liu \cite{BL2014}. Wang, Yang and Zheng~\cite{WYZ2026} determined the asymptotic minimum
for every number of edges between $\mathrm{ex}(n,K_{r+1})+1$ and
$\mathrm{ex}(n,K_{r+2})$.

\smallskip 
We study this jump at a spectral threshold that involves only the
number of edges. Let $\lambda(G)$ denote the spectral radius of $G$,
that is, the largest eigenvalue of its adjacency matrix.
Nosal~\cite{nosal} proved that every triangle-free graph with $m$
edges satisfies $\lambda(G)\le\sqrt m$. We refer to $\sqrt m$ as the
\emph{Nosal threshold}. Following Li, Liu and Zhang~\cite{LLZ2026JCTB}, we call
$\lambda(G)>\sqrt m$ the \emph{Nosal condition} and a graph that
satisfies it a {\bf \color{blue} \emph{Nosal graph}}. 
For the convenience of readers, we postpone the numerous known results on Nosal graphs to Section \ref{subsec-1.2}. 
 Since every graph satisfies
$\lambda(G)\ge2m/n$, the Mantel condition
$e(G)>\lfloor n^2/4\rfloor$ implies the Nosal condition. The converse
is false, because the Nosal condition constrains neither the order
nor the degree distribution. For instance, every complete bipartite
graph $K_{a,b}$ satisfies $\lambda^2=ab=m$, however unbalanced its
parts are. The boundary of the Nosal condition therefore contains
graphs with far fewer than $n^2/4$ edges, and~\eqref{e-1} yields no
sharp bound in terms of $m$. These bipartite graphs are triangle-free,
so $f$ vanishes on them, just as it vanishes on $T_{n,2}$. The
spectral form of the Erd\H{o}s--Tuza problem asks how large $f$ must
become once a graph crosses this threshold. We study it through the
function
\[
 \psi(m) :=\min\bigl\{f(G): e(G)=m,\ G\text{ is }K_4\text{-free},\
 \lambda(G)>\sqrt m\bigr\}.
\] 
A single book gives a first lower bound. We call the common
neighborhood $N(u)\cap N(v)$ of an edge $uv$ its \emph{page}, so that
the \emph{booksize} $\bk(G)$, the largest number of triangles sharing an
edge, is the largest size of a page. All pairs inside a page of a
$K_4$-free graph are saturating non-edges, so
$f(G)\ge\binom{\bk(G)}{2}$. The sharp bound $\bk(G)\ge \frac{1}{3}\lambda(G)$ of
Zhao, You, Zeng and Zhang~\cite{ZYZZ2026} then yields
$\psi(m)\ge \frac{1}{18}m-O(\sqrt m\,)$. This argument uses one page, whereas the
Balogh--Liu proof \cite{BL2014} uses all three pages of a triangle, which are
pairwise disjoint. In the spectral setting, the eigenvalue equations
involve both the sizes of these pages and the weights of the Perron
vector on them. We control the two quantities at once, which closes
the gap between $\frac{1}{18}m $ and the true value $\frac{1}{9}m$.

\subsection{Main results}\label{subsec-1.1}
Let $q(G)$ be the number of vertices of $G$ that lie in at least one
triangle. Throughout this paper, we put $\tF(G)=f(G)+ \frac 12 q(G)$. Our first result is a sharp
inequality for this corrected count.
\begin{theorem}\label{thm-1.1}
Let $G$ be a $K_4$-free graph with $m$ edges. If $\lambda(G)>\sqrt m$,
then
\[
 f(G)+\frac{1}2 q(G) \ > \ \frac{1}{9} \lambda(G)^2.
\]
The constant $1/9$ is best possible.
\end{theorem}

The correction $q(G)/2$ makes the count scale exactly like
$\lambda(G)^2$. Let $G^{(s)}$ be the graph obtained from $G$ by
replacing each vertex by an independent set of size $s$ and each edge
by a complete bipartite graph $K_{s,s}$. Saturating pairs in distinct classes
contribute $s^2f(G)$, and a pair inside one class is saturating
exactly when the corresponding vertex of $G$ lies in a triangle. Hence, we have 
$f(G^{(s)})=s^2f(G)+q(G)\binom s2$, which gives
$\tF(G^{(s)})=s^2\tF(G)$, while
$\lambda(G^{(s)})=s\lambda(G)$.

\smallskip 
Theorem~\ref{thm-1.1} has the following unconditional
form, which bounds $\lambda(G)$ for $K_4$-free graphs. 

\begin{corollary}\label{cor-1.2}
Every $K_4$-free graph $G$ with $m$ edges satisfies
\[
 \lambda(G)^2\ \le \ \max\bigl\{m,\,9\tF(G)\bigr\}.
\]
\end{corollary}
\noindent
Indeed, if $\lambda(G)^2>m$, then Theorem~\ref{thm-1.1} gives
$9\tF(G)>\lambda(G)^2$. For a triangle-free graph, $f=q=0$, so the
maximum equals $m$ and Corollary~\ref{cor-1.2} reduces to Nosal's
inequality. In this sense, Theorem~\ref{thm-1.1} extends Nosal's
theorem from triangle-free graphs to $K_4$-free graphs.

\smallskip 
The corrected count also bounds $f(G)$ itself, with a linear error term.

\begin{theorem}\label{thm-1.3}
Every $K_4$-free Nosal graph $G$ with $m$ edges satisfies
\[
 f(G) \ > \ \frac{1}{9} \lambda(G)^2 -\lambda(G),\qquad
 f(G) \ >\ \frac m9-\sqrt m.
\]
Moreover, for every $m\ge500$,
\[
 \frac m9-\sqrt m \ < \ \psi(m) \ < \  \frac m9+\sqrt m .
\]
In particular, we have $\lim_{m\to \infty}\psi(m) / m= 1/9$.
\end{theorem}

Theorem~\ref{thm-1.3} determines the size of the jump at the Nosal
threshold, and all
its constants are explicit.
Over $K_4$-free graphs with $m$ edges, the minimum of $f$ is $0$ on
the boundary $\lambda(G)=\sqrt m$, where the star $K_{1,m}$ attains
it, and it jumps to $( \frac 19+o(1))m$ just above the boundary.

\smallskip 
For large spectral radius, we sharpen the linear error term as follows.

\begin{theorem}\label{thm-1.4}
Every $K_4$-free Nosal graph $G$ with sufficiently large spectral
radius satisfies
\[
 f(G) \ > \ \frac{1}{9} \lambda(G)^2 -\frac{1}3 \lambda(G).
\]
\end{theorem}

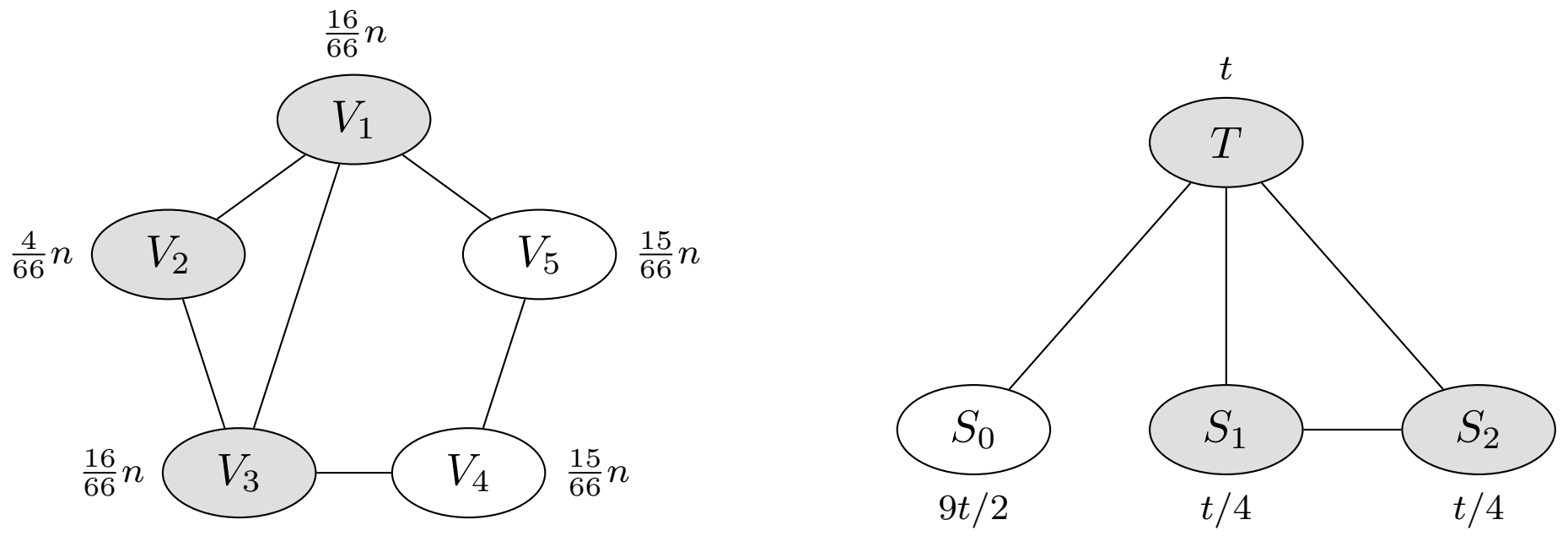
\begin{figure}[ht]
\centering
\begin{tikzpicture}[scale=0.9,
 cls/.style={draw,ellipse,minimum width=12mm,minimum height=7mm,
             inner sep=1pt,font=\small},
 sat/.style={cls,fill=gray!25},
 lab/.style={font=\footnotesize},
 sz/.style={font=\scriptsize}]
\begin{scope}
\node[sat] (v1) at (90:1.7)  {$V_1$};
\node[sat] (v2) at (162:1.7) {$V_2$};
\node[sat] (v3) at (234:1.7) {$V_3$};
\node[cls] (v4) at (306:1.7) {$V_4$};
\node[cls] (v5) at (18:1.7)  {$V_5$};
\draw (v1)--(v2)--(v3)--(v4)--(v5)--(v1);
\draw (v1)--(v3);
\node[sz,above=0pt of v1] {$\frac{16}{66}n$};
\node[sz,left=0pt of v2] {$\frac{4}{66}n$};
\node[sz,left=0pt of v3] {$\frac{16}{66}n$};
\node[sz,right=0pt of v4] {$\frac{15}{66}n$};
\node[sz,right=0pt of v5] {$\frac{15}{66}n$};
\node[lab,align=center] at (0,-2.75)
 {(a) The Balogh--Liu graph:\\ $e(G)>n^2/4$,\quad
  $f(G)=(\frac2{33}+o(1))n^2$};
\end{scope}
\begin{scope}[xshift=7.6cm]
\node[sat] (T)  at (0,1.5)  {$T$};
\node[cls] (S0) at (-2.2,-1)  {$S_0$};
\node[sat] (S1) at (0,-1)   {$S_1$};
\node[sat] (S2) at (2.2,-1)   {$S_2$};
\draw (T)--(S0); \draw (T)--(S1); \draw (T)--(S2); \draw (S1)--(S2);
\node[sz,above=0pt of T] {$t$};
\node[sz,below=0pt of S0] {$9t/2$};
\node[sz,below=0pt of S1] {$t/4$};
\node[sz,below=0pt of S2] {$t/4$};
\node[lab,align=center] at (0,-2.75)
 {(b) $G_{t,1/4}$: $n=6t$,\quad $m=\lambda^2=81t^2/16$,\\
  $\lambda=9t/4$,\quad $f=\lambda^2/9-\lambda/3$};
\end{scope}
\end{tikzpicture}
\caption{Two extremal constructions in the classical and spectral settings. 
Each node is an independent set
with the size shown, up to additive constants in (a), and a line
stands for a complete bipartite graph. In both graphs, the saturating
non-edges are exactly the pairs inside the shaded classes. For
general $p$, the graph $G_{t,p}$ has $|T|=t$, $|S_0|=(4+2p)t$ and
$|S_1|=|S_2|=pt$.}
\label{fig-1}
\end{figure}

\noindent 
{\bf The extremal graphs.} 
The graphs $G_{t,p}$ consist of an independent
set $T$ of size $t$ joined completely to a set $S$ of size
$(4+4p)t$. Inside $S$, a complete bipartite graph joins two classes
of size $pt$, and the other vertices of $S$ have no neighbors in $S$.
Each of these graphs satisfies $\lambda^2=m$, and the choice $p=1/4$
minimizes $\tF/\lambda^2$. For $4\mid t$, the graph $G_{t,1/4}$ in
Figure~\ref{fig-1}(b) satisfies $\tF=\lambda^2/9$ and
$f=\lambda^2/9-\lambda/3$ exactly. Deleting one vertex of $S$ with no
neighbors in $S$ gives a Nosal graph, and the resulting graphs satisfy
$\tF/\lambda^2\to1/9$; see Section~\ref{sec-5} for detailed computations.

\smallskip 
These graphs differ from the extremal graphs of Balogh and Liu \cite{BL2014}. For
$m=\lfloor n^2/4\rfloor+1$, the bound~\eqref{e-1} reads
$f(G)\ge(8/33-o(1))m$, which is more than twice $m/9$. The graphs that
reach $1/9$ satisfy neither the Mantel condition nor the
\emph{spectral Mantel condition} $\lambda(G)>\lambda(T_{n,2})$.
Indeed, $G_{t,1/4}$ has $6t$ vertices, $81t^2/16<9t^2$ edges and
spectral radius $9t/4<3t=\lambda(T_{6t,2})$.

\smallskip
These results show that the jump at the Nosal threshold
is a new phenomenon and not a rescaled form of the
Balogh--Liu bound (\ref{e-1}). We highlight three features.
First, we determine the size of the jump:
$\psi(m)=\frac19m+O(\sqrt m\,)$ with explicit constants,
which doubles the bound $\frac1{18}m$ from a single largest book.
Second, Theorem~\ref{thm-1.1} extends Nosal's theorem with a
sharp constant, and the count $\tF$ scales exactly like
the power $\lambda^2$ under blow-ups. Third, the pages of a triangle
in $G_{t,1/4}$ have sizes in the ratio $4:1:1$; see Figure~\ref{fig-1}(b). The three
pages therefore play different roles, and a sharp proof must
keep all of them at once. In Section~\ref{subsec-1.3}, we shall explain how our proof keeps all
three pages, and why the ratio $4:1:1$ survives the averaging.

\smallskip 
Theorem~\ref{thm-1.4} is close to optimal. The boundary
graphs $G_{t,1/4}$ as defined above satisfy $\lambda^2=m$ and
$f=\frac19\lambda^2-\frac13\lambda$. Among Nosal graphs, the
graphs $G'$ of Section~\ref{subsec-5.2} satisfy 
$f<\frac19\lambda^2-\frac{23}{81}\lambda $. Thus, 
the best coefficient of the linear term in Theorem~\ref{thm-1.4} lies between
$\frac{23}{81}\approx0.284$ and $\frac13$.

\subsection{Related work}\label{subsec-1.2}

Saturating non-edges underlie the well-known saturation problem of Erd\H{o}s,
Hajnal and Moon~\cite{EHM1964} and the weak saturation problem
of Bollob\'as~\cite{Boll1967} and Lov\'{a}sz \cite{Lov1977}. 
For a fixed graph $F$, a graph $G$ is called 
\emph{$F$-saturated} if it is $F$-free and adding any non-edge of $G$ creates
a copy of $F$. 
Erd\H{o}s, Hajnal and Moon \cite{EHM1964} showed
that every $K_t$-saturated graph on $n$ vertices has at least $\binom n2-\binom{n-t+2}2$ edges, and this bound is attained only at the join $K_{t-2}\vee \overline{K}_{n-t+2}$.  
More generally, the saturation number of $F$ is defined to be the minimum number of edges
of an $F$-saturated graph on $n$ vertices. K\'aszonyi and
Tuza~\cite{KT1986} proved that this number is at most linear in
$n$ for every fixed $F$. Kritschgau, Methuku, Tait and
Timmons~\cite{KMTT2020} studied the minimum number of copies of a
fixed graph in an $F$-saturated graph. 
We refer to the comprehensive survey~\cite{CFFS2021}
for known results.

\smallskip 
The spectral version of the Erd\H{o}s--Hajnal--Moon problem asks for  
 the minimum spectral radius of a $K_t$-saturated graph on $n$ vertices, 
 and whether $K_{t-2}\vee \overline{K}_{n-t+2}$ is the unique extremal graph.  
Kim, Kim, Kostochka and O~\cite{KKKO2020} first studied the case $t=3$.  
Kim, Kostochka, O, Shi and Wang~\cite{KKOSW2023} settled the case $t=4$, 
and they proposed a sufficient criterion based on $2$-walks for every $t\ge 3$.  
By verifying this criterion, Wang and Hou~\cite{wanghou} settled the cases $t=5,6$, 
and Ai, Liu, O and Zhang~\cite{ALOZ2025} independently settled the case $t=5$. 
Chen and Zhang \cite{CZ2026} solved the remaining cases. 
 We mention that apart from cliques, 
 there are few results on spectral saturation problems for other forbidden graphs. 
 Incidentally, the spectral problem of this paper is quite different: we do not require $G$ to
be $K_4$-saturated, and we count its $K_4$-saturating non-edges
instead. 

\smallskip 
Spectral extremal problems in terms of the number $m$ of edges go back to Brualdi and
Hoffman~\cite{BH1985} and Stanley~\cite{stanley}, who bounded
the spectral radius of a graph $G$ with $m$ edges.
Rowlinson~\cite{rowlinson} later determined the extremal graphs. 
Nikiforov~\cite{Nik2002} extended Nosal's theorem to
$\lambda(G)^2\le (1- \frac{1}{r})2m$ for every $K_{r+1}$-free graph $G$, and  
 later he \cite{Nik2006,Nik2009} showed that 
equality holds if and only if $G$ is either a complete bipartite
graph for $r=2$ or a regular complete $r$-partite graph for $r\ge 3$, 
possibly with isolated vertices. 
Zhai, Lin and Shu~\cite{ZLS2021} studied such problems
systematically under the name of Brualdi--Hoffman--Tur\'an problems; 
see \cite{LZS2024,LLZ-edge-spectral,LLZ-edge-color-critical} for recent related results. 
In addition, the papers~\cite{LNW2021,zhaishu,LFP2024,LLZ2025}
refine Nosal's bound for nonbipartite graphs. 

\smallskip 
Much recent work quantifies how much additional structure
a Nosal graph must contain beyond a single triangle. In terms of the order, spectral
supersaturation goes back to Bollob\'as and
Nikiforov~\cite{BN2007}, who bounded the number of cliques
from below in terms of $\lambda(G)$ and $n$. At
the Nosal threshold, Ning and Zhai~\cite{ningzhai} proved that every Nosal graph 
contains at least $\lfloor \frac{1}{2}(\sqrt{m} -1) \rfloor$ triangles. 
They \cite{ningzhai2} also showed that every Nosal graph 
with $m\ge3.6\times10^9$ edges contains at least $\frac{1}{2000}m^2$ copies of $C_4$. 
Li, Liu and Zhang~\cite{LLZ2026JCTB} determined the optimal constant $\frac{1}{8}$.  
Zhai, Lin and Shu~\cite{ZLS2021} proved that every Nosal graph contains a complete bipartite graph $K_{2,t}$ with $t>\frac{1}{4} \sqrt m$.
Li, Liu and Zhang~\cite{LLZ-edge-spectral} improved this to $t\ge \frac{1}{2}\sqrt m-O(1)$, and random graphs show that the constant $\frac{1}{2}$ is best possible. 
 For books, Nikiforov~\cite{Nik2021} proved
that every Nosal graph satisfies $\bk (G) > \frac{1}{12}m^{1/4}$. 
By a probabilistic argument, 
Li, Liu and Zhang~\cite{LLZ2026JCTB} first obtained the correct order by showing that
$\bk(G)> \frac{1}{24} \sqrt m$. Zhai, Li and Lou~\cite{zhaililou} improved the constant to $\frac{1}{9}$. 
Recently, Zhao, You, Zeng and
Zhang~\cite{ZYZZ2026} proved the sharp inequality
$\bk(G)\ge \frac{1}{3}\lambda(G) > \frac{1}{3}\sqrt{m}$.
They also showed that every Nosal graph contains at least $\lambda(G)$ edges in
triangles, which confirms a conjecture of Li, Feng and
Peng~\cite{LFP2025}. 
The book $B_2$ is also called the \emph{kite} and denoted by $C_4^+$, since it arises from $C_4$ by adding a chord.
Chen, Li and Tang~\cite{CLT2026} showed that every Nosal graph contains at least $\big(\frac{1}{8} -o(1) \big)m$ kites, and the constant $\frac{1}{8}$ is best possible.  
Table~\ref{tab-1} compares these bounds with Theorem~\ref{thm-1.3}.

\vspace{-1mm}
\begin{table}[H]
\centering
\renewcommand{\arraystretch}{1.2}
\setlength{\tabcolsep}{10pt}
\begin{tabular}{lcc}
\toprule
Forced substructure  & Above Mantel: $e(G)> \lfloor n^2/4\rfloor$ & Above Nosal: $\lambda (G) > \sqrt{m}$ \\
\midrule
 $\#C_3$ & $\lfloor \frac{n}{2} \rfloor$    \cite{Erd1955,Erdos1964}
 &  $\big\lfloor \frac{1}{2}(\sqrt{m} -1) \big\rfloor$
 \cite{ningzhai}  \\
 $\#C_4$ & $\Omega (n^4) $  & $\big(\frac{1}{8} - o(1) \big)m^2 $ \cite{LLZ2026JCTB} \\
 Largest $K_{2,t}$ & $\Omega(n)$ & $ \frac{1}{2}\sqrt{m} -O(1)$ \cite{LLZ-edge-spectral}  \\
 Booksize $\bk(G)$ & $\lceil \frac{n}{6} \rceil $  \cite{KN1979}
 &  $ \frac{1}{3}\sqrt{m}$ \cite{ZYZZ2026}  \\
 $\#$Triangular edges & $2\lfloor \frac{n}{2}\rfloor +1$  \cite{EFR1992} & $ \sqrt{m} $ \cite{ZYZZ2026} \\
 Largest chordal subgraph & $\frac{3}{2}n-1$  \cite{EGOZ1989}  & $\Omega (\sqrt{m}\,)$ \cite{LLZ2026JCTB} \\
 $\#C_4^+$ & $\binom{\lfloor n/2\rfloor}{2}$
 \cite{Mub2010} &  $\big(\frac{1}{8} -o(1) \big) m$ \cite{CLT2026} \\
 $\# B_t$ for fixed $t\ge 1$ & $\binom{\lfloor n/2\rfloor}{t}$
 \cite{Mub2010} & $\big(\frac{1}{t!2^t} \!-\! o(1) \big)m^{t/2}$ \cite{CLT2026} \\
\midrule
 $\#K_4$-saturating non-edges & $\frac{2}{33}n^2 \! +\! \Theta (n)$
 \cite{BL2014} & $\frac 19 m \! + \! \Theta(\sqrt m\,)$ (Theorem~\ref{thm-1.3}) \\
\bottomrule
\end{tabular}
\caption{Lower bounds forced above the Mantel and Nosal thresholds.
The last row is the new result and assumes that $G$ is $K_4$-free.
 Here $C_4^+$ is the $4$-cycle with one chord, and $B_t$ is the
book of size $t$. A triangular edge is an edge that lies in a
triangle, and the chordal row counts edges.}
 \label{tab-1}
\end{table}

We would like to point out that all earlier known results on Nosal graphs, as stated in Table \ref{tab-1}, 
bound substructures that the graph contains. In contrast, the function $\psi (m)$ of this paper counts pairs
that the graph misses, and its extremal structures are different from those of~\eqref{e-1}.

\subsection{Proof overview and organization}\label{subsec-1.3}

\noindent
\textbf{Proof overview.} The proof has to overcome three difficulties. 
First, the Nosal
condition does not involve the order. The graphs $G_{t,1/4}$
have only $\frac9{64}n^2$ edges, so~\eqref{e-1} does not
apply to them directly. Second,  our lower bounds for $f$
come from the sizes of pages, whereas the eigenvalue
equations involve the Perron weights on these pages. 
We must control the two quantities at once. Third, and most importantly, the decisive
estimates leave no room for loss. The Nosal condition gives
only $\lambda^2-m>0$, and these estimates must be tight on
$G_{t,1/4}$, where every triangle has pages in the ratio
$4:1:1$. In particular, we cannot replace the three pages of
a triangle by the largest one.

\smallskip
The constant $1/9$ comes from all three pages of a triangle.
We first explain this for a join $G=I_a\vee H$, where $I_a$
is an independent set on $a$ vertices and $H$ is
triangle-free. Scale the Perron vector $\by$ of $G$ to equal
$1$ on $I_a$, and choose an edge $uv$ of $H$ that maximizes
$y_u+y_v$. Together with a vertex of $I_a$, this edge forms
a triangle whose pages are disjoint and have sizes $a$,
$b=d_H(u)$ and $c=d_H(v)$. Proposition~\ref{prop-2.1} shows
that
\[
 \lambda<2a+\frac{b+c}{2},\qquad
 \tF(G)\ge\frac{a^2+b^2+c^2}{2}.
\]
Since $(2a+\frac{b+c}2)^2\le\frac92(a^2+b^2+c^2)$ by the
Cauchy--Schwarz inequality, it follows that
$\lambda^2<9\tF(G)$. Equality in the Cauchy--Schwarz step
holds exactly when $a:b:c=4:1:1$. The maximizing edge has
two properties. First, its surplus $y_u+y_v-1$ is positive,
because the surpluses of all edges of $H$ sum to
$\lambda^2-e(G)>0$. This is where the Nosal condition
enters. Second, every neighbor of $u$ in $H$ has coordinate
at most $y_v$, and every neighbor of $v$ in $H$ has
coordinate at most $y_u$.

In a general graph, a vertex $u^*$ of largest Perron
coordinate plays the role of a vertex of $I_a$, and $H$
becomes $G[N(u^*)]$. The vertices outside $N(u^*)$ need not
be adjacent to all of $N(u^*)$. The Perron weight that a
vertex of $N(u^*)$ misses outside $N[u^*]$ lowers the
surplus of its edges. After this correction, the Nosal
condition gives a positive surplus only on average, and the
maximizing edge need not have it. Thus, the two properties
may hold on different edges.

\smallskip
To get around this difficulty, we compare two families of
numbers on all edges of $H$ instead of choosing one edge.
The numbers $s_{uv}$ encode the second property. Double
counting gives $\sum_{uv\in E(H)}s_{uv}=0$, and
$s_{uv}\le0$ at a maximizing edge. So this identity is an
averaged form of maximality, and it needs no choice of edge.
The numbers $g_{uv}$ encode the corrected surplus, and the
Nosal condition gives $\sum_{uv\in E(H)}g_{uv}>0$
(Lemma~\ref{lem-3.1}). The essential observation in our
proof is that if the claimed bound fails, then
$s_{uv}\ge(r/2)g_{uv}$ on every edge $uv$, where $r:= (1+y(W))/\lambda $ is the ratio
of the Perron weight outside $N(u^*)$ to the weight inside
$N(u^*)$. The coefficient $r/2$ is the same on all edges, so
summing over $E(H)$ gives
$0\ge\frac r2\sum_{uv\in E(H)}g_{uv}>0$, which leads to a contradiction.
Each edge enters with the sizes of its own pages, so the
ratio $4:1:1$ survives the averaging. On the graph $G_{t,1/4}$, where
$y_u=y_v=1/2$ on every edge of $H$, both $s_{uv}$ and
$g_{uv}$ vanish, so the comparison holds with equality
there and leaves no room for loss. 

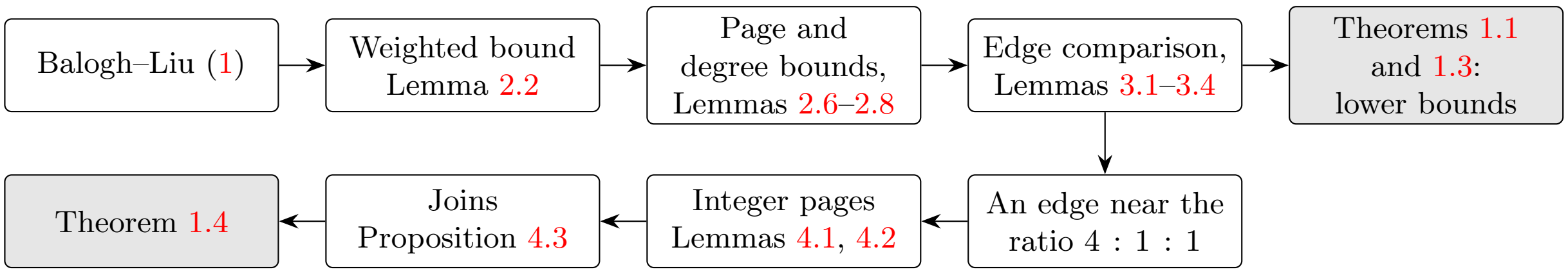
\begin{figure}[ht]
\centering
\begin{tikzpicture}[
 box/.style={draw,rounded corners=2pt,align=center,font=\footnotesize,
             inner sep=3pt,minimum height=9.5mm,text width=26mm},
 goal/.style={box,fill=gray!20},
 arr/.style={-{Stealth[length=2mm]}}]
\node[box] (a1) at (0,0) {Balogh--Liu~\eqref{e-1}};
\node[box] (a2) at (3.3,0) {Weighted bound\\ Lemma~\ref{lem-2.2}};
\node[box] (a3) at (6.6,0) {Page and degree bounds, Lemmas~\ref{lem-2.6}--\ref{lem-2.8}};
\node[box] (a4) at (9.9,0) {Edge comparison, Lemmas~\ref{lem-3.1}--\ref{lem-3.4}};
\node[goal] (a5) at (13.2,0) {Theorems~\ref{thm-1.1} and~\ref{thm-1.3}: lower bounds};
\node[box] (b4) at (9.9,-1.6) {An edge near the ratio $4:1:1$};
\node[box] (b3) at (6.6,-1.6) {Integer pages\\ Lemmas~\ref{lem-4.1},~\ref{lem-4.2}};
\node[box] (b2) at (3.3,-1.6) {Joins\\ Proposition~\ref{prop-4.3}};
\node[goal] (b1) at (0,-1.6) {Theorem~\ref{thm-1.4}};
\draw[arr] (a1)--(a2);\draw[arr] (a2)--(a3);\draw[arr] (a3)--(a4);
\draw[arr] (a4)--(a5);\draw[arr] (a4)--(b4);\draw[arr] (b4)--(b3);
\draw[arr] (b3)--(b2);\draw[arr] (b2)--(b1);
\end{tikzpicture}
\caption{The structure of the proofs.}
\label{fig-2}
\end{figure}

\smallskip
The proof of the lower bounds has three steps; see
Figure~\ref{fig-2}. First, a blow-up weighted by the squared
Perron coordinates transfers~\eqref{e-1} to the spectral
setting (Lemma~\ref{lem-2.2}). This is where the Balogh--Liu theorem enters. 
Second, an exact identity for the spectral surplus $\lambda^2-m$ localizes
the Nosal condition at $u^*$ (Lemma~\ref{lem-2.4}). If the
claimed bound fails, then the weighted bound and the page
counts restrict $r$ and the degrees in $H$
(Lemmas~\ref{lem-2.6} -- \ref{lem-2.8}). Third,
Lemmas~\ref{lem-3.2} -- \ref{lem-3.4} reduce the comparison on
each edge to a quadratic polynomial in $y_u$ and $y_v$, and
show that this polynomial is nonnegative. 
The two inputs play different roles. Lemma~\ref{lem-2.2} is
not tight on $G_{t,1/4}$, where its two sides differ by a
factor of about $2$, and we use it for rough estimates such
as $r<7/10$. For $r\le1/3$, the comparison holds
automatically. Near the extremal value $r=4/9$, the page
counts of Lemma~\ref{lem-2.6} decide the comparison. Thus, 
the sharp constant $1/9$ comes from the pages, while the
Balogh--Liu theorem supplies rough bounds on $r$ and on the
degrees in $H$.

\smallskip
Our method contains two parts that extend beyond the count $f(G)$.
The first is a weighted blow-up according to the Perron coordinates, adapted from the
 probabilistic blow-up of Li, Liu and
 Zhang~\cite{LLZ2026JCTB}. A lower bound above
the Mantel threshold, for a count with a quadratic leading
term under blow-ups, becomes a Perron-weighted bound above
the Nosal threshold (Lemma~\ref{lem-2.2}).
Section~\ref{sec-6} gives the version for larger cliques.
The second is a zero-sum comparison, which replaces an
extremal choice of edge. It combines a family of edge
numbers that has zero sum and is nonpositive at a
maximizing edge with a family that has positive sum by the
spectral condition. We compare the two families on each edge
with one coefficient, which the extremal graph determines.
We expect this scheme to apply to other counts that depend
on several pages of a triangle.

\smallskip
Finally, the similar comparison drives the proof of
Theorem~\ref{thm-1.4}. If $f(G)\le\lambda^2/9-\lambda/3$,
then the comparison fails at some edge of $H$. The pages
$A,B,C$ of the triangle formed by this edge and $u^*$ have
sizes close to the ratio $4:1:1$ and satisfy
$2|A|+(|B|+|C|)/2>\lambda$. Integrality then fixes the page
sizes and shows that every saturating non-edge lies in these
pages (Lemma~\ref{lem-4.1}). The eigenvalue equations force
$G$ to be a join, and we treat joins directly
(Proposition~\ref{prop-4.3}).

\medskip
\noindent
\textbf{Organization.} 
The rest of the paper is organized as follows. 
Section~\ref{sec-2} develops the three-page
counts, the transference and the estimates for Perron coordinates. Section~\ref{sec-3}
proves the lower bounds in Theorems~\ref{thm-1.1}
and~\ref{thm-1.3} by edge averaging, and Section~\ref{sec-4} proves
Theorem~\ref{thm-1.4}. Section~\ref{sec-5} gives the extremal family
and the upper bound in Theorem~\ref{thm-1.3}, and Section~\ref{sec-6}
discusses further problems for the interested reader and provides possible extensions of our
method. 

\section{Three pages and spectral localization}\label{sec-2}

To begin with, we introduce some notation needed later. For a graph
$G$, we write $V(G)$ and
$E(G)$ for its vertex set and edge set, $N_G(v)$ for the neighborhood
of a vertex $v$, and $N_G[v]=N_G(v)\cup\{v\}$. We omit the subscript
when the graph is clear from the context. For $X\subseteq V(G)$, let
$G[X]$ be the induced subgraph, and put $e(X)=e(G[X])$,
$N_X(v)=N(v)\cap X$ and $d_X(v)=|N_X(v)|$. Let $Q(G)$ be the set of
vertices of $G$ that lie in triangles, so that $q(G)=|Q(G)|$ in the
notation of Section~\ref{subsec-1.1}.

\subsection{Why the constant is \texorpdfstring{$1/9$}{1/9}}
\label{subsec-2.1}

Let $I_a$ denote an independent set of $a\ge1$ vertices, and let
$\vee$ denote the graph join. The next proposition proves
Theorem~\ref{thm-1.1} for joins $I_a\vee H$, where the comparison of
three pages by the Cauchy--Schwarz inequality appears in its simplest form.

\begin{proposition}\label{prop-2.1}
Let $G=I_a\vee H$, where $H$ is triangle-free, and suppose that
$\lambda=\lambda(G)>\sqrt{e(G)}$. Then $H$ has an edge $uv$ such that
$b=d_H(u)$ and $c=d_H(v)$ satisfy
\[
 \lambda<2a+\frac{b+c}{2},\qquad
 \tF(G)\ge\frac{a^2+b^2+c^2}{2}.
\]
In particular, $\lambda^2<9\tF(G)$.
\end{proposition}
\begin{proof}
If $H$ had no edge, then $G$ would be complete bipartite and
$\lambda^2=e(G)$. So $H$ has an edge, and $G$ is connected. The
vertices of $I_a$ have the same neighborhood, so we may normalize the
positive eigenvector $\by=(y_v)$ of $\lambda$ to equal $1$ on $I_a$.
Then $\sum_{v\in V(H)}y_v=\lambda$, and summing the eigenvalue
equations over $V(H)$ yields
$\lambda^2=a|V(H)|+\sum_{uv\in E(H)}(y_u+y_v)$. Since
$e(G)=a|V(H)|+e(H)$, this identity becomes
\[
 \sum_{uv\in E(H)}(y_u+y_v-1)=\lambda^2-e(G)>0.
\]
Choose an edge $uv$ of $H$ that maximizes $s=y_u+y_v$, so that $s>1$,
and put $x=y_u$, $z=y_v$, $b=d_H(u)$ and $c=d_H(v)$. By maximality,
every neighbor of $u$ in $H$ has coordinate at most $z$, and every
neighbor of $v$ in $H$ has coordinate at most $x$. Hence
$\lambda x\le a+bz$ and $\lambda z\le a+cx$. Multiply these
inequalities by $x$ and $z$, and add them. Since
$x^2+z^2\ge s^2/2>s/2$ and $xz\le(x^2+z^2)/2$, this leads to
\[
 \lambda(x^2+z^2)\le as+(b+c)xz
 <\left(2a+\frac{b+c}{2}\right)(x^2+z^2).
\]
For any $w\in I_a$, the sets $I_a$, $N_H(u)$ and $N_H(v)$ are the
common neighborhoods of the three edges of the triangle $wuv$. They
are independent and pairwise disjoint, because $H$ is triangle-free.
Every pair inside one of them is saturating, and each of their
vertices lies in a triangle. A set of size $k$ among them therefore
contributes $\binom k2+k/2=k^2/2$ to $\tF(G)$, which proves the second
inequality. Finally, by the Cauchy--Schwarz inequality,
\[
 \lambda^2<\left(2a+\frac{b+c}{2}\right)^2 
 \le\frac92(a^2+b^2+c^2)\le9\tF(G). \qedhere
\]
\end{proof}

The Cauchy--Schwarz step above is tight exactly when
$a:b:c=4:1:1$, which is the ratio in the family of
Section~\ref{sec-5}. The general proof keeps the three
page counts, and a zero-sum comparison over all edges of a
neighborhood replaces the choice of a maximizing edge. 

\subsection{Transference from the order to the size}
\label{subsec-2.2}

For a connected graph $G$, the Perron--Frobenius theorem provides an
eigenvector for $\lambda(G)$ with all positive entries. We call it the
\emph{Perron vector} of $G$ and its entries the Perron coordinates. For the
rest of this section, let $G$ be a connected $K_4$-free Nosal graph,
and write $m=e(G)$ and $\lambda=\lambda(G)$. Let $\by=(y_v)_{v\in V(G)}$ be its
Perron vector, normalized by $\|\by\|_{\infty} = \max_{v\in V(G)} y_v =1$. 
We write $\|\by\|_1=\sum_{v\in V(G)} y_v$ and $\|\by\|_2 =(\sum_{v\in V(G)} y_v^2)^{1/2}$ throughout
the paper.

\smallskip 
The next lemma transfers the Balogh--Liu bound~\eqref{e-1} to the spectral setting in terms of its Perron vector. The main
ingredient is a blow-up that uses the squared Perron coordinates as
vertex weights. Under the Nosal condition, this blow-up has more
edges than the Mantel threshold for its order. 
A probabilistic blow-up argument according to the Perron coordinates was used by Li, Liu and Zhang \cite[Lemma 2.6]{LLZ2026JCTB} for bounding the booksize of Nosal graphs from below.

\begin{lemma}\label{lem-2.2}
The Perron vector $\by$ satisfies
\[
 \sum_{uv\in\Sat(G)}y_u^2y_v^2+
 \frac12\sum_{v\in Q(G)}y_v^4
 \ge\frac2{33}\|\by\|_2^4 .
\]
\end{lemma}
\begin{proof}
Put $z_v:=y_v^2/\|\by\|_2^2$ for every $v\in V(G)$, so these positive weights sum to $1$. By
the Cauchy--Schwarz inequality, the eigenvalue equation and the Nosal condition, we have
\[
 p:=\sum_{uv\in E(G)}z_uz_v
 \ge\frac{1}{m\|\by\|_2^4} \biggl(\sum_{uv\in E(G)}y_uy_v\biggr)^2
 =\frac{\lambda^2}{4m}>\frac14.
\]
We first treat positive rational weights $z_v$ with sum $1$ and
$p>1/4$. Let $s$ tend to infinity through multiples of a common
denominator. Replace each vertex $v$ by an independent set of size
$sz_v$ and each edge by the corresponding complete bipartite graph.
The resulting graph $G_s$ is still $K_4$-free on $s$ vertices with
$ps^2>s^2/4$ edges. Delete edges to reach exactly
$\lfloor s^2/4\rfloor+1$ edges, and apply~\eqref{e-1}. Every saturating
non-edge of this spanning subgraph is already a saturating non-edge of
$G_s$: if the pair were an edge of $G_s$, it would form a $K_4$ with
its witness edge. Thus $f(G_s)\ge(2/33-o(1))s^2$.

Saturating pairs in distinct classes correspond to those of $G$, and a
pair inside the class of $v$ is saturating exactly when $v\in Q(G)$.
Hence
\[
 \frac{f(G_s)}{s^2}
 =\sum_{uv\in\Sat(G)}z_uz_v+
   \frac12\sum_{v\in Q(G)}z_v^2
   -\frac1{2s}\sum_{v\in Q(G)}z_v.
\]
Letting $s\to\infty$, we obtain
$\sum_{\Sat(G)}z_uz_v+\frac12\sum_{Q(G)}z_v^2\ge2/33$ for every rational vector $\bm{z}$.
Rational approximation preserves the strict inequality $p>1/4$, so by continuity the same
bound holds for every positive weight vector with sum $1$ and
$p>1/4$. Multiplying by $\|\by\|_2^4$ proves the lemma.
\end{proof}

The above proof uses two features of the saturation count: it has a
quadratic leading term under blow-ups, and Balogh--Liu's result \eqref{e-1} bounds it from
below above the Mantel threshold.

\subsection{Localization at a Perron vertex}\label{subsec-2.3}

We now localize the problem at a vertex of largest Perron coordinate.
Keep $G$, $\lambda$ and $\by$ as above. For $X\subseteq V(G)$, write
$y(X)=\sum_{v\in X}y_v$. Choose a vertex $u^*$ with $y_{u^*}=1$ and
put
\[
 U=N(u^*),\qquad W=V(G)\setminus N[u^*],\qquad H=G[U],
 \qquad r=\frac{1+y(W)}{\lambda}.
\]
The graph $H$ is triangle-free, since $G$ is $K_4$-free. The parameter $r$ records the normalized
Perron weight outside $U$. On the extremal graph $G_{t,1/4}$ in Figure~\ref{fig-1}(b), it
equals $4/9$.

Recall that the page of an edge $uv$ is defined as $N(u)\cap N(v)$. By the three
pages of a triangle, we mean the pages of its three edges.

\begin{lemma}\label{lem-2.3}
In a $K_4$-free graph, every page is independent, and every pair of
distinct vertices in a page is a saturating non-edge. The three pages
of a triangle are pairwise disjoint, and their vertices lie in
triangles. Consequently, if these pages have sizes $p_1,p_2,p_3$, then
\[
 f(G)\ge\binom{p_1}2+\binom{p_2}2+\binom{p_3}2,\qquad
 \tF(G)\ge\frac{p_1^2+p_2^2+p_3^2}{2}.
\]
\end{lemma}
\begin{proof}
An edge inside $N(u)\cap N(v)$ would form a $K_4$ with $u,v$, and
adding an edge between two vertices of this page creates a $K_4$ with
$u,v$. A vertex in two pages of a triangle would be adjacent to all
three vertices of that triangle, again creating a $K_4$. The two
bounds follow, since a page of size $k$ contributes
$\binom k2+k/2=k^2/2$ to $\tF(G)$.
\end{proof}

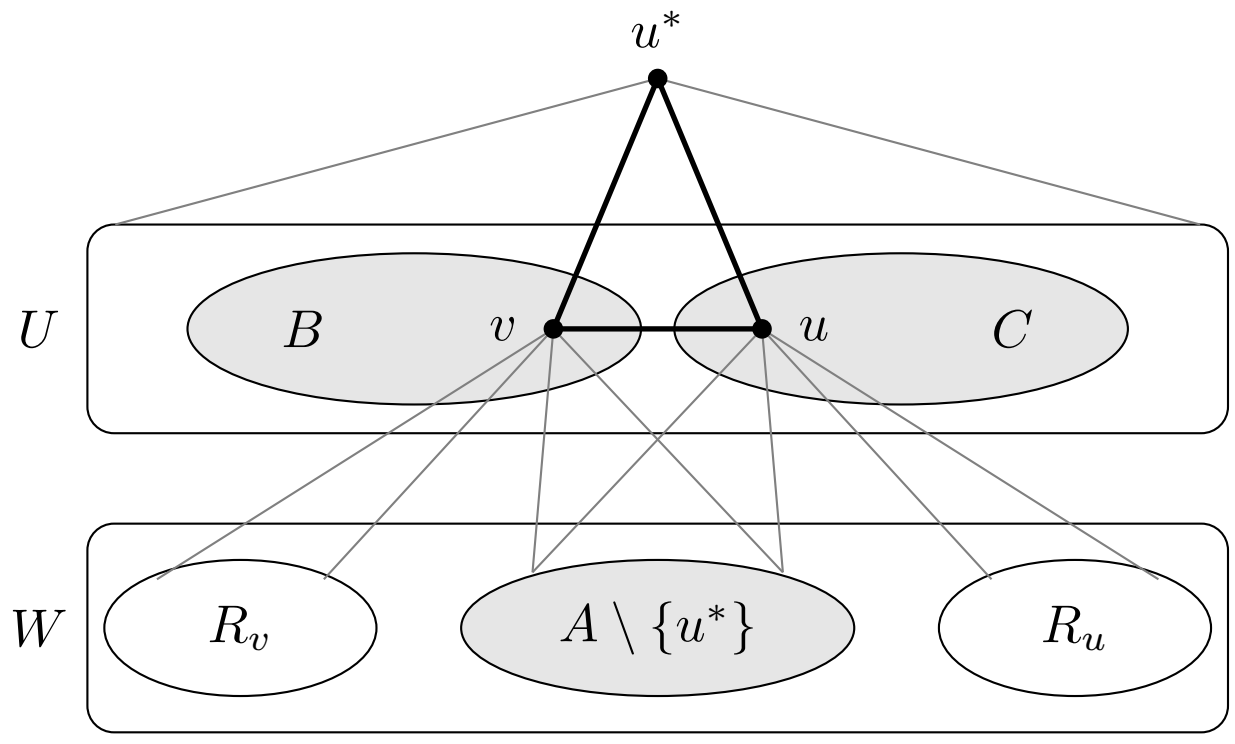
\begin{figure}[ht]
\centering
\begin{tikzpicture}[scale=0.92,
 vtx/.style={circle,fill=black,inner sep=1.3pt},
 pg/.style={draw,fill=gray!20,ellipse,inner sep=0pt},
 ex/.style={draw,ellipse,inner sep=0pt},
 fan/.style={gray,line width=.4pt},
 tag/.style={font=\footnotesize}]
\draw[rounded corners=5pt] (-4.1,-.05) rectangle (4.1,1.45);
\draw[rounded corners=5pt] (-4.1,-2.2) rectangle (4.1,-.7);
\node[font=\small] at (-4.45,.7) {$U$};
\node[font=\small] at (-4.45,-1.45) {$W$};
\coordinate (s) at (0,2.5);
\draw[fan] (s)--(-3.9,1.45); \draw[fan] (s)--(3.9,1.45);
\node[pg,minimum width=30mm,minimum height=10mm] (B) at (-1.75,.7) {};
\node[pg,minimum width=30mm,minimum height=10mm] (C) at (1.75,.7) {};
\node[font=\small] at (-2.55,.7) {$B$};
\node[font=\small] at (2.55,.7) {$C$};
\node[pg,minimum width=26mm,minimum height=9mm] (A) at (0,-1.45) {};
\node[ex,minimum width=18mm,minimum height=9mm] (Rv) at (-3.0,-1.45) {};
\node[ex,minimum width=18mm,minimum height=9mm] (Ru) at (3.0,-1.45) {};
\node[font=\small] at (0,-1.45) {$A\setminus\{u^*\}$};
\node[font=\small] at (-3.0,-1.45) {$R_v$};
\node[font=\small] at (3.0,-1.45) {$R_u$};
\coordinate (v) at (-.75,.7);
\coordinate (u) at (.75,.7);
\draw[fan] (v)--(-3.6,-1.1); \draw[fan] (v)--(-2.4,-1.1);
\draw[fan] (v)--(-.9,-1.05); \draw[fan] (v)--(.9,-1.05);
\draw[fan] (u)--(3.6,-1.1); \draw[fan] (u)--(2.4,-1.1);
\draw[fan] (u)--(.9,-1.05); \draw[fan] (u)--(-.9,-1.05);
\draw[line width=1pt] (s)--(v)--(u)--cycle;
\node[vtx,label={[font=\small]above:$u^*$}] at (s) {};
\node[vtx,label={[font=\small,xshift=-1pt]left:$v$}] at (v) {};
\node[vtx,label={[font=\small,xshift=1pt]right:$u$}] at (u) {};
\end{tikzpicture}
\caption{The sets at an edge $uv$ of $H=G[U]$. The pages of the
triangle $u^*uv$ are shaded: $A$ contains $u^*$, $B$ contains $v$, and
$C$ contains $u$. They are independent and pairwise disjoint, and
every pair inside a page is a $K_4$-saturating non-edge. Note that on the graph $G_{t,1/4}$, the
page sizes have ratio $4:1:1$, and the sets $R_u$ and $R_v$ are empty.}
\label{fig-3}
\end{figure}

For an edge $uv\in E(H)$, the three pages of the triangle $u^*uv$ are
\[
 A=N(u)\cap N(v),\qquad B=N_U(u),\qquad C=N_U(v).
\]
Here $B,C\subseteq U$ and $A\subseteq W\cup\{u^*\}$. Put
$R_u=N_W(u)\setminus N(v)$ and $R_v=N_W(v)\setminus N(u)$. Then $N(u)$
is the disjoint union of $A,B,R_u$, and $N(v)$ is the disjoint union
of $A,C,R_v$; see Figure~\ref{fig-3}.

The next identity expresses the spectral surplus
$\lambda^2-m$ through the edges of $H=G[U]$ and the vertices of $W$.
Chen, Li and Tang~\cite{CLT2026} used this identity to
bound the booksize, and Zhao, You, Zeng and
Zhang~\cite{ZYZZ2026} used it to count the triangular
edges of Nosal graphs. An earlier variant of this
identity appeared in the work of Zhai, Li and
Lou~\cite{zhaililou}. 

\begin{lemma}\label{lem-2.4}
The set $U$ satisfies $y(U)=\lambda$, and
\[
 \lambda^2-m=\mu-\nu>0,\] 
 where 
\[ \mu :=\sum_{uv\in E(H)}(y_u+y_v-1) \quad \text{and} \quad 
 \nu:=\sum_{w\in W}d_U(w)(1-y_w)+e(W).
\]
In particular, $H$ has an edge.
\end{lemma}
\begin{proof}
The equation at $u^*$ reads $\lambda=y(U)$. Summing the equations at
the vertices of $U$, we get
\[
 \lambda^2=|U|+\sum_{uv\in E(H)}(y_u+y_v)
                 +\sum_{w\in W}d_U(w)y_w.
\]
Now subtract $m=|U|+e(H)+\sum_{w\in W}d_U(w)+e(W)$. Finally, we have $\mu>0$
because $\nu\ge0$.
\end{proof}

\begin{lemma}\label{lem-2.5}
The number $r$ satisfies $0<r<1$ and $\|\by\|_2^2>2r\lambda$. Moreover,
for every edge $uv\in E(H)$ with pages $A,B,C$ defined as above,
\[
 y(A)+y(B)+y(C)\ge\lambda(y_u+y_v-r).
\]
\end{lemma}
\begin{proof}
Every coordinate of $\bm{y}$ is at most $1$, so $d(v)\ge\lambda y_v$ for every
vertex $v$. Summing the eigenvalue equations, we find
$\lambda\|\by\|_1=\sum_vd(v)y_v\le2m<2\lambda^2$. Since
$\|\by\|_1=\lambda(1+r)$, it follows that $r<1$. The sum
$\sum_v(d(v)-\lambda y_v)(1-y_v)$ is nonnegative, and expanding it
shows that
\[
 \lambda\|\by\|_2^2\ge2\lambda\|\by\|_1-2m
 >2\lambda^2(1+r)-2\lambda^2=2r\lambda^2.
\]

For the second assertion, the sets $A\setminus\{u^*\},R_u,R_v$ are
disjoint subsets of $W$, so $y(A)+y(R_u)+y(R_v)\le r\lambda$. Now, 
the eigenvalue equations at $u$ and $v$ give  
$\lambda y_u = y(B)+ y(A) + y(R_u)$ and $\lambda y_v = y(C) + y(A) + y(R_v)$. 
Adding the two equations and using the inequality $y(A)+y(R_u)+y(R_v)\le r\lambda$, 
we obtain the second assertion of the lemma. 
\end{proof}

\subsection{Page and degree bounds}\label{subsec-2.4}

The next three lemmas form the harder part of this
section. They show
that a failure of our bounds forces small pages, a small Perron weight
$r$ and small degrees in $H$. Throughout this subsection, we assume
that one of the following holds:
\[
 \text{\rm(i)}\ \ \tF(G)\le\frac{\lambda^2}{9};\qquad
 \text{\rm(ii)}\ \ f(G)+\lambda\le\frac{\lambda^2}{9};\qquad
 \text{\rm(iii)}\ \ f(G)+\frac\lambda3\le\frac{\lambda^2}{9}
 \ \text{ and }\ \lambda\ge500.
\]
Cases (i) and (ii) serve Theorems~\ref{thm-1.1} and~\ref{thm-1.3}, and
case (iii) serves Theorem~\ref{thm-1.4}. Put $\tau=0$ in cases (i) and
(ii), and $\tau=1/\lambda$ in case (iii). For an edge $uv\in E(H)$
with pages $A,B,C$ of $u^*uv$, put $b=|B|/\lambda$ and
$c=|C|/\lambda$. Finally, let $\alpha=\max_{v\in U}y_v$.

\begin{lemma}\label{lem-2.6}
The pages of every triangle of $G$ have sizes $p_1,p_2,p_3$ with
$p_1^2+p_2^2+p_3^2\le(2/9+\tau/3)\lambda^2$. In particular, for every
edge $uv\in E(H)$, we have $b<1/2$ and $b+c<3/4$.
\end{lemma}
\begin{proof}
In case (i), this is Lemma~\ref{lem-2.3}. In cases (ii) and (iii),
$f(G)\le\lambda^2/9$ and $\lambda\ge9$. Let $k=p_1+p_2+p_3$. By
Lemma~\ref{lem-2.3} and convexity, $k(k-3)/6\le f(G)$. If $k>\lambda$,
then $\lambda(\lambda-3)/6<\lambda^2/9$, that is, $\lambda<9$. So
$k\le\lambda$, and Lemma~\ref{lem-2.3} implies
\[
 \frac{p_1^2+p_2^2+p_3^2}{2}\le f(G)+\frac k2
 \le f(G)+\frac\lambda2\le\frac{\lambda^2}{9}+\frac{\tau\lambda^2}{6}.
\]
The last two claims follow from $\tau\le1/500$.
\end{proof}

\begin{lemma}\label{lem-2.7}
For every edge $uv\in E(H)$,
\[
 \frac8{33}r^2+\frac{1-\alpha^4}{2}(b^2+c^2)
 <\frac19+\frac{2\tau}3.
\]
In particular, $r<7/10$.
\end{lemma}
\begin{proof}
The proof relies on the following observation: in Lemma~\ref{lem-2.2},
the pairs and vertices inside $B$ and $C$ have weight at most
$\alpha^4$ rather than $1$. Let $Y$ denote the left side of
Lemma~\ref{lem-2.2}. By Lemmas~\ref{lem-2.2} and~\ref{lem-2.5},
we have $Y>(8/33)r^2\lambda^2$. Therefore, the lemma follows from the bound
\[
 Y\le\Bigl(\frac19+\frac{2\tau}3\Bigr)\lambda^2
    -\frac{1-\alpha^4}{2}\lambda^2(b^2+c^2).
\]
In case (i), this is immediate from the observation, which yields
$Y\le\tF(G)-(1-\alpha^4)(|B|^2+|C|^2)/2$. In cases (ii) and (iii), we
treat the vertex terms separately. Since $y_v^4\le\alpha^3y_v$ on $U$
and $y_v^4\le y_v$ on all other vertices, we have
\[
 \sum_{v\in Q(G)}y_v^4\le\alpha^3y(U)+y(V(G)\setminus U)
 =\lambda(\alpha^3+r).
\]
Together with the observation for pairs, this implies
\[
 Y\le f(G)-\frac{1-\alpha^4}{2}\lambda^2(b^2+c^2)
  +\frac\lambda2\bigl(\alpha^3+r+(1-\alpha^4)(b+c)\bigr).
\]
By the weighted AM--GM inequality, we get $\alpha^3\le(1+3\alpha^4)/4$, so $b+c<3/4$ implies
$\alpha^3+(1-\alpha^4)(b+c)\le1$. With $r<1$, the last term is less
than $\lambda$. In both cases, we get 
$f(G)+\lambda\le\lambda^2/9+2\tau\lambda^2/3$, and the bound on $Y$
follows. The final claim holds because
$r^2<(33/8)(1/9+2\tau/3)<49/100$.
\end{proof}

\begin{lemma}\label{lem-2.8}
The coordinates on $U$ satisfy $\alpha^2<5/8$. Consequently, for every
edge $uv\in E(H)$,
\[
 r^2+\frac54(b^2+c^2)<\frac{11}{24}+3\tau.
\]
\end{lemma}
\begin{proof}
A vertex of $U$ with a large coordinate would force $r$ to be large,
and Lemma~\ref{lem-2.7} forbids this. To make this precise, suppose
that $\alpha^2\ge5/8$, and choose $u\in U$ with $y_u=\alpha$. Since
$r<7/10<\alpha$, the equation at $u$ yields
$y(N_U(u))\ge\lambda(\alpha-r)>0$, so $u$ has a neighbor in $H$. Let
$v$ be its neighbor in $H$ of largest coordinate, put $\beta=y_v$, and
let $A,B,C$ be the pages of $u^*uv$. The same equation implies
$\alpha-r\le b\beta$. Coordinates on $A,B,C$ are at most
$1,\beta,\alpha$, respectively, and $\beta\le\alpha$. By
Lemmas~\ref{lem-2.5} and~\ref{lem-2.6}, and by the Cauchy--Schwarz inequality, we obtain 
\[
 \alpha+\beta-r
 \le\sqrt{\Bigl(\frac29+\frac\tau3\Bigr)(1+2\alpha^2)}
 \le\alpha\sqrt{\frac45+\frac{6\tau}5}<\frac9{10}\alpha,
\]
where the middle inequality uses $1/\alpha^2\le8/5$. Thus
$\beta<r-\alpha/10$. Inserting this bound into $\alpha-r\le b\beta$ and
solving for $r$, we get
\[
 r>\alpha\frac{1+b/10}{1+b}
 \ge\alpha\Bigl(\frac9{10}-\frac{2b}{5}\Bigr),
\]
where the difference of the last two terms is
$2\alpha (b-1/2)^2/ \big(5(1+b) \big) \ge0$. Now apply Lemma~\ref{lem-2.7} to this edge,
discard $c^2$, and use this lower bound on $r$. The result is
$1/9+2\tau/3>t$, where
\[
 t=\frac8{33}\alpha^2\Bigl(\frac9{10}-\frac{2b}{5}\Bigr)^2
 +\frac{1-\alpha^4}{2}b^2.
\]
The number $t$ is concave in $\alpha^2$, so its minimum on $[5/8,1]$
occurs at an endpoint. At $\alpha^2=1$, we find
$t>(8/33)(7/10)^2>9/80$, since $b<1/2$. At $\alpha^2=5/8$, the
inequalities $5/33>3/20$ and $39/128>3/10$ show that
\[
 t\ge\frac3{20}\Bigl(\frac9{10}-\frac{2b}{5}\Bigr)^2+\frac3{10}b^2
 =\frac9{80}+\frac{81}{250}\Bigl(b-\frac16\Bigr)^2.
\]
Thus $t\ge9/80$ in all cases, whereas $1/9+2\tau/3<9/80$ because
$\tau\le1/500$. This contradiction proves that $\alpha^2<5/8$, which is the
first assertion of the lemma.

For the second assertion, multiply the inequality of
Lemma~\ref{lem-2.7} by $33/8$ and use $1-\alpha^4>39/64$, together
with the numerical fact $(33/16)(39/64)=1287/1024>5/4$.
\end{proof}

\section{Edge averaging}\label{sec-3}

In this section, we prove the lower bounds in Theorems~\ref{thm-1.1}
and~\ref{thm-1.3}. The argument is based on a comparison of two
families of numbers on the edges of $H$. One family sums to zero by
double counting, and the other has positive sum by the Nosal
condition (Lemma~\ref{lem-3.1}). If a bound fails, then the first
family dominates $r/2$ times the second on every edge, which is
impossible. The delicate part of the proof lies in this edge
inequality: Lemma~\ref{lem-3.2} reduces it to a quadratic polynomial,
and Lemmas~\ref{lem-3.3} and~\ref{lem-3.4} show that this polynomial is
nonnegative.

Let $G$ be a connected $K_4$-free Nosal graph with $m$ edges. 
In the sequel, we keep the previous notation
of Section~\ref{subsec-2.3}. For each vertex $v\in U$, put
\[
 b_v=\frac{d_H(v)}{\lambda},\qquad
 t_v=\frac{y(N_H(v))}{\lambda},\qquad
 \omega_v=\frac{y(W\setminus N(v))}{\lambda}.
\]
The number $\omega_v$ records the Perron weight in $W$ that is missing
from the neighborhood of $v$, and the equation at $v$ gives
$t_v=y_v-r+\omega_v$. For each edge $uv$ of $H$, put
\[
 s_{uv}=y_ut_u+y_vt_v-(b_u+b_v)y_uy_v,\qquad
 g_{uv}=2(y_u+y_v-1)-\frac{y_u\omega_u}{b_u}
        -\frac{y_v\omega_v}{b_v}.
\]
\begin{lemma}\label{lem-3.1}
The two families satisfy
\[
 \sum_{uv\in E(H)}s_{uv}=0 \quad \text{and} \quad  \sum_{uv\in E(H)}g_{uv}>0.
\]
\end{lemma}
\begin{proof}
Counting at the endpoints of edges, we see that
$\sum_{uv\in E(H)}(y_ut_u+y_vt_v)=\sum_{v\in U}d_H(v)y_vt_v$. In
$\sum_{uv\in E(H)}(b_u+b_v)y_uy_v$, a vertex $u$ contributes
$b_uy_u\,y(N_H(u))=d_H(u)y_ut_u$. So the two parts of the first sum
cancel, which proves the first assertion.

For the second assertion, we interchange sums and use the eigenvalue
equations on $W$:
\[
 \lambda\sum_{v\in U}y_v\omega_v
 =\sum_{w\in W}y_w\bigl(\lambda-y(N_U(w))\bigr)
 =\sum_{w\in W}y(N_U(w))(1-y_w)
   +\sum_{ww'\in E(G[W])}(y_w+y_{w'}).
\]
Since $y(N_U(w))\le d_U(w)$ and $y_w+y_{w'}\le2$, the right side is at
most $2\nu$. A vertex $v$ of positive degree in $H$ contributes
$-d_H(v)y_v\omega_v/b_v=-\lambda y_v\omega_v$ to
$\sum_{uv\in E(H)}g_{uv}$. Hence, by Lemma~\ref{lem-2.4},
\[
 \sum_{uv\in E(H)}g_{uv}
 \ge2\mu-\lambda\sum_{v\in U}y_v\omega_v\ge2(\mu-\nu)>0.\qedhere
\]
\end{proof}

We next reduce the inequality $s_{uv}\ge(r/2)g_{uv}$ to the
nonnegativity of a quadratic polynomial in the coordinates of $u$ and
$v$.
\begin{lemma}\label{lem-3.2}
The endpoints of each edge of $H$ can be named $u,v$ so that the
following holds. Put $a=y(N(u)\cap N(v))/\lambda$, $b=b_u$, $c=b_v$,
$h=b+c$, $x=y_u$, $z=y_v$ and
\[
 w=(r-a)\Bigl(1+\frac{r}{2b}\Bigr).
\]
Then $a\le r$ and
\[
 s_{uv}-\frac r2g_{uv}\ge\rho(x,z),\quad\text{where}\quad
 \rho(x,z)=x^2+z^2-hxz-2r(x+z)+r+wx.
\]
\end{lemma}
\begin{proof}
Let $A=N(u)\cap N(v)$. The inclusion $A\subseteq W\cup\{u^*\}$ gives
$y(A)\le1+y(W)=r\lambda$, so $a\le r$. Moreover,
$W\setminus A$ is the union of $W\setminus N(u)$ and
$W\setminus N(v)$, so
$\omega_u+\omega_v\ge y(W\setminus A)/\lambda=r-a$. Substituting
$t_u=x-r+\omega_u$ and $t_v=z-r+\omega_v$, we obtain
\[
 s_{uv}-\frac r2g_{uv}
 =x^2+z^2-hxz-2r(x+z)+r
  +x\omega_u\Bigl(1+\frac r{2b}\Bigr)
  +z\omega_v\Bigl(1+\frac r{2c}\Bigr).
\]
The sum of the last two terms is at least
$(r-a)\min\{x(1+r/(2b)),\,z(1+r/(2c))\}$. Name the endpoints so that
the first entry is the minimum.
\end{proof}

\begin{lemma}\label{lem-3.3}
Let $0\le h<2$, $0<r\le1-h^2/4$ and $w\ge0$, and let $\rho$ be as in
Lemma~\ref{lem-3.2}. Define
\[
 k=4r+h-2,\qquad\gamma(w)=(2+h)r(2w-k)-w^2.
\]
If $\rho(x,z)<0$ for some $x,z\ge0$, then $w<2r$, $k>0$ and
$\gamma(w)<0$.
\end{lemma}
\begin{proof}
If $w\ge2r$, then for $x,z\ge0$,
\[
 \rho(x,z)\ge x^2-hxz+z^2-2rz+r
 \ge\Bigl(1-\frac{h^2}4\Bigr)z^2-2rz+r
 \ge r-\frac{r^2}{1-h^2/4}\ge0.
\]
So $w<2r$. Completing the square in $x+z$ and $z-x$, we find
\[
 \rho(x,z)=\frac{2-h}{4}\Bigl(x+z-\frac{4r-w}{2-h}\Bigr)^2
 +\frac{2+h}{4}\Bigl(z-x-\frac{w}{2+h}\Bigr)^2
 +\frac{\gamma(w)}{4-h^2}.
\]
Thus $\rho(x,z)<0$ forces $\gamma(w)<0$. This in turn forces $k>0$,
since $k\le0$ and $w<2r$ would imply
$\gamma(w)\ge w\bigl(2(2+h)r-w\bigr)\ge0$.
\end{proof}

\begin{lemma}\label{lem-3.4}
Let $a,c\ge0$ and $b,r>0$ satisfy $a\le r$ and
\[
 a^2+b^2+c^2\le\frac29,\qquad
 r^2+\frac54(b^2+c^2)\le\frac{11}{24}.
\]
Put $h=b+c$ and $w=(r-a)(1+r/(2b))$, and let $k,\gamma$ be as in
Lemma~\ref{lem-3.3}. If $k>0$ and $w\le2r$, then $\gamma(w)>0$ holds
for this value of $w$.
\end{lemma}
\begin{proof}
The trick of the proof is to compare $w$ with the fixed multiple
$11k/18$ of $k$. Since $\gamma$ is concave, it suffices to establish
\[
 \text{(a)}\ \ \gamma(2r)>0,\qquad \quad 
 \text{(b)}\ \ \gamma\bigl(\tfrac{11k}{18}\bigr)>0,\qquad \quad 
 \text{(c)}\ \ \tfrac{11k}{18}\le w.
\]
Indeed, $w$ then lies between $11k/18$ and $2r$. Parts (a) and (b) are
short computations, and part (c), where the hypothesis
$a^2+b^2+c^2\le2/9$ enters, requires the most work. Throughout, we use
three simple facts: $h^2\le2(b^2+c^2)\le4/9$, $r<7/10$, and
$r>(2-h)/4\ge1/3$, where the last of these bounds restates the
assumption $k>0$ of the lemma.

\emph{Proof of {\rm(a)}.} By the definition, we see that 
$\gamma(2r)=r(4-h^2-4r)$, and this is positive.

\emph{Proof of {\rm(b)}.} Direct expansion shows that
\[
 \gamma\Bigl(\frac{11k}{18}\Bigr)
 =\frac{k}{324}\bigl(242-340r-121h+72hr\bigr).
\]
So we need $340r+121h-72hr<242$. Suppose first that $r\ge1/2$. Then
$72hr\ge36h$. Since $r^2+(5/8)h^2\le r^2+(5/4)(b^2+c^2)\le11/24$,
the Cauchy--Schwarz inequality yields
\[
 340r+85h\le\sqrt{\Bigl(340^2+\frac85\cdot85^2\Bigr)\frac{11}{24}}
 =\sqrt{\frac{127160\cdot11}{24}}<242.
\]
If $r<1/2$, the left side, which increases with $r$, is less than
$170+85h\le170+85\cdot\frac23<242$.

\emph{Proof of {\rm(c)}.} Put $\sigma=2-h-4a$, so that
$k=4(r-a)-\sigma$. Expanding, we obtain
\[
 \frac29\sigma=\Bigl(\frac29-a^2-b^2-c^2\Bigr)
 +\Bigl(a-\frac49\Bigr)^2+\Bigl(b-\frac19\Bigr)^2
 +\Bigl(c-\frac19\Bigr)^2.
\]
Hence $\sigma\ge0$, and more precisely
\[
 \sigma\ge\frac92\Bigl(a-\frac49\Bigr)^2
          +\frac92\Bigl(b-\frac19\Bigr)^2.
\]
Substituting $k=4(r-a)-\sigma$ into $w-11k/18$, we arrive at
\[
 2b\Bigl(w-\frac{11k}{18}\Bigr)
 =(r-a)\Bigl(r-\frac{26b}{9}\Bigr)+\frac{11}{9}b\sigma.
\]
This is nonnegative if $r\ge26b/9$. Suppose that $r<26b/9$, so that
$b>9r/26>3/26$. By AM--GM,
$(r-a)(26b/9-r)\le(26b/9-a)^2/4$, so the claim reduces to
$(44/9)b\sigma\ge(26b/9-a)^2$. Put $a'=4-9a$ and $b'=9b-1$. By the
lower bound for $\sigma$, it is enough to prove that
\[
 198(b'+1)(a'^2+b'^2)\ge t^2,\qquad\text{where}\quad
 t=26b'+9a'-10=81\Bigl(\frac{26b}9-a\Bigr)>0.
\]
If $26b'\le10$, then $0<t\le9a'$, and the inequality holds because
$b'+1=9b>1$. If $26b'>10$, then the Cauchy--Schwarz inequality yields
\[
 t^2=\Bigl(\bigl(26-\tfrac{10}{b'}\bigr)b'+9a'\Bigr)^2
 \le\Bigl(\bigl(26-\tfrac{10}{b'}\bigr)^2+81\Bigr)(a'^2+b'^2).
\]
So we need $(26b'-10)^2+81b'^2\le198b'^2(b'+1)$. The difference of the
two sides is
\[
 198b'^3-559b'^2+520b'-100=(b'-1)^2(198b'-163)+63-4b'.
\]
Here $5/13<b'<13/4$, because $b\le\sqrt2/3$. Thus $63-4b'>50$. The
first term is nonnegative for $b'\ge163/198$, and otherwise it exceeds
$-(8/13)^2\cdot87>-34$. So the difference is positive.
\end{proof}

As promised, we now deduce the lower bounds in our main results.

\begin{proof}[Proof of the lower bounds in Theorems~\ref{thm-1.1}
and~\ref{thm-1.3}]
Suppose first that $G$ is connected and that one of the two bounds
fails, so that case (i) or (ii) of Section~\ref{subsec-2.4} holds with
$\tau=0$. Consider an edge of $H$, with endpoints named as in
Lemma~\ref{lem-3.2}, and let $A,B,C$ be the pages of $u^*uv$. Then
$b=|B|/\lambda$, $c=|C|/\lambda$ and $a\le|A|/\lambda$. So
Lemmas~\ref{lem-2.6} and~\ref{lem-2.8} show that $a,b,c,r$ satisfy the
hypotheses of Lemma~\ref{lem-3.4}. In particular, $h^2\le4/9$ and
$r<7/10$, so $r\le1-h^2/4$. By Lemmas~\ref{lem-3.3} and~\ref{lem-3.4},
$\rho(x,z)\ge0$, and Lemma~\ref{lem-3.2} then yields
$s_{uv}\ge(r/2)g_{uv}$. Since $r>0$ does not depend on the edge,
summing over $E(H)$ contradicts Lemma~\ref{lem-3.1}.

For a disconnected graph, choose a component $G_0$ of spectral radius
$\lambda$. Since $\lambda^2>e(G)\ge e(G_0)$, the connected result
applies to $G_0$. As $f$ and $q$ are additive over components, both
bounds extend from the component $G_0$ to the whole graph $G$.

Finally, we deduce the bound in terms of $m$. If $m<81$, then
$m/9-\sqrt m<0$. For $m\ge81$, the function $x^2/9-x$ increases on
$[\sqrt m,\infty)$. Since $\lambda>\sqrt m$, it follows that
$f(G)>\lambda^2/9-\lambda>m/9-\sqrt m$. This proves the lower bound on
$\psi(m)$ as well. Section~\ref{subsec-5.2} shows that the constant
$1/9$ in Theorem~\ref{thm-1.1} is best possible.
\end{proof}

\section{The linear term at large spectral radius}\label{sec-4}

Write $\phi(t)=t^2/9-t/3$. Our main tools in proving
Theorem~\ref{thm-1.4} are the edge comparison of Section~\ref{sec-3}
and the integrality of page sizes. We begin with two combinatorial
lemmas. We then treat joins $I_a\vee H$ with $H$ triangle-free.

\begin{lemma}\label{lem-4.1}
Let $a,b,c$ be nonnegative integers, and put
\[
 \ell=2a+\frac{b+c}{2},\qquad
 f_0=\binom a2+\binom b2+\binom c2,\qquad
 j=a-2(b+c).
\]
Then $f_0-\phi(\ell)=j(j+3)/18+(b-c)^2/4$. Moreover, suppose that
$7\le\lambda<\ell$ and that an integer $f$ satisfies
$f_0\le f\le\phi(\lambda)$. Then
\[
 b=c=d\ge2,\qquad a\in\{4d-2,4d-1\},\qquad f=f_0,
\]
and $\varepsilon=\ell-\lambda$ satisfies
$\varepsilon(\ell+\lambda-3)\le1$.
\end{lemma}
\begin{proof}
The identity follows by substituting $a=2(b+c)+j$. Its first summand
is at least $-1/9$, with equality for $j=-1,-2$ alone, while the second
is either zero or at least $1/4$. Since $\phi$ increases on
$[7,\infty)$, we have $f_0\le\phi(\lambda)<\phi(\ell)$. This forces
$b=c=d$ and $j\in\{-1,-2\}$, so $\phi(\ell)-f_0=1/9$. By integrality,
$f=f_0$, and the bound $\phi(\ell)-\phi(\lambda)\le1/9$ is the claimed
inequality for $\varepsilon$. Finally, we have $d\ge2$, because $d\le1$ would
imply $\ell=2a+d\le9d-2\le7$.
\end{proof}

\begin{lemma}\label{lem-4.2}
Let $G$ be a $K_4$-free graph, let $u^*\in V(G)$, and put
$H=G[N(u^*)]$. Let $uv\in E(H)$, and let $A,B,C$ be the pages of
$u^*uv$, with $B=N_H(u)$ and $C=N_H(v)$. If every saturating non-edge
of $G$ lies inside $A$, $B$ or $C$, then $H$ is the disjoint union of
a bipartite graph with parts $B,C$, a matching and a set of isolated
vertices.
\end{lemma}
\begin{proof}
The sets $B,C$ are disjoint and independent, and $u\in C$, $v\in B$.
Suppose that $w\in V(H)\setminus(B\cup C)$ has a neighbor $w'\in B$.
Then $w,u$ are nonadjacent, and their common neighbors include the
adjacent vertices $u^*,w'$. So $wu$ is a saturating non-edge outside
the three pages, a contradiction. The same argument applies to $C$.
Thus, no edge of $H$ joins $B\cup C$ to its complement.

Now let $x$ be a vertex of $H-(B\cup C)$ with two neighbors there.
These neighbors are nonadjacent, because $H$ is triangle-free. They
form a saturating non-edge witnessed by $x$ and $u^*$, and this pair
lies outside the three pages, because $A$ is disjoint from $V(H)$.
\end{proof}

For joins, we next show that the linear term holds as soon as $\lambda\ge7$.
\begin{proposition}\label{prop-4.3}
Let $G=I_a\vee H$, where $a\ge1$ and $H$ is triangle-free. Suppose that
$\lambda=\lambda(G)>\sqrt{e(G)}$ and $\lambda\ge7$. Then the linear bound
$f(G)>\phi(\lambda)$ holds for $G$.
\end{proposition}
\begin{proof}
Suppose that $f(G)\le\phi(\lambda)$. Our first goal is to show that
$H=K_{d,d}\cup sK_1$ for some integers $d,s$. To this end, take the
edge $uv\in E(H)$ of Proposition~\ref{prop-2.1}, with $b=d_H(u)$,
$c=d_H(v)$ and $\lambda<\ell=2a+(b+c)/2$. For $u^*\in I_a$, the pages
of the triangle $u^*uv$ are $I_a$, $B=N_H(u)$ and $C=N_H(v)$, so 
$f(G)\ge\binom a2+\binom b2+\binom c2$ by Lemma~\ref{lem-2.3}. Now, 
Lemma~\ref{lem-4.1} implies $b=c=d\ge2$, $a\in\{4d-2,4d-1\}$ and
equality in this bound. Thus, every saturating non-edge lies inside
$I_a$, $B$ or $C$, and Lemma~\ref{lem-4.2} yields
$H=H_0\cup pK_2\cup sK_1$, where $H_0$ is bipartite with parts $B,C$.
In particular, we have $\Delta(H)\le d$, $e(H)\le d^2+p$, $a\ge3d\ge6$ and
$\ell=2a+d$. Moreover, we see that $\varepsilon=\ell-\lambda>0$ satisfies
$\varepsilon(\ell+\lambda-3)\le1$, so $\varepsilon<1$,
$\lambda>\ell-1>d$ and $\lambda \varepsilon <1$.

Next we measure how far the degrees in $H$ are from $d$. Normalize the
Perron vector to equal $1$ on $I_a$. As in
Proposition~\ref{prop-2.1}, we see that  
$\mu=\sum_{ij\in E(H)}(y_i+y_j-1)=\lambda^2-e(G)>0$. The equation at a
vertex of largest coordinate in $H$ shows that this coordinate is at
most $\beta=a/(\lambda-d)$. With $k_i=d_H(i)$, the equation at $i$
then implies
\[
 y_i\le\frac{a+k_i\beta}{\lambda}=\beta-\frac{(d-k_i)\beta}{\lambda}.
\]
Summing over edge endpoints, we get
$0<\mu\le(2\beta-1)e(H)-(\beta/\lambda)z$, where
$z=\sum_ik_i(d-k_i)\ge0$. Since $2\beta-1=\beta\varepsilon/a$, this
means that
\[
 z<\frac{\lambda\varepsilon}{a}e(H)<\frac{e(H)}{a}.
\]
If $p\ge1$, the matching contributes $2p(d-1)\ge pd$ to $z$, whereas
$e(H)/a\le(d^2+p)/(3d)<pd$, a contradiction. So $p=0$. If
$H_0\ne K_{d,d}$, then a nonadjacent pair in $B\times C$ provides two
vertices of degree between $1$ and $d-1$. Then $z\ge2(d-1)\ge d$,
while $e(H)/a\le d/3$, again a contradiction. This proves that
$H=K_{d,d}\cup sK_1$, as desired.

The final step is an integrality argument. Now $e(G)=a(2d+s)+d^2$, so
$\ell^2-e(G)=a(4a+2d-s)$ is an integer multiple of $a$. It is positive,
because $e(G)<\lambda^2<\ell^2$. The Perron vector is constant on the
orbits of the automorphism group, so the eigenvalue equations reduce
to $\lambda=2ad/(\lambda-d)+as/\lambda$. Multiplying by $\lambda$ and
using $2a+d-\lambda=\varepsilon$, we find
$\lambda^2-e(G)=d^2\varepsilon/(\lambda-d)$. Since
$\lambda-d>\ell-1-d>a\ge d$ and $\ell<\lambda+1$, this leads to
\[
 a\le\ell^2-e(G)
 =\varepsilon\Bigl(\ell+\lambda+\frac{d^2}{\lambda-d}\Bigr)
 <3\lambda\varepsilon<3,
\]
which is contrary to $a\ge6$. This completes the proof. 
\end{proof}

We are now in a position to prove the last of our main results.

\begin{proof}[Proof of Theorem~\ref{thm-1.4}]
To start the proof, we first outline the main steps. The aim of our
proof is to show that a counterexample must be a join, which
Proposition~\ref{prop-4.3} excludes. First of all, we establish that a
counterexample has an edge of $H$ whose pages have ratio close to
$4:1:1$. The key point is that, by Lemma~\ref{lem-3.4}, the comparison
of Section~\ref{sec-3} can fail only at such an edge. Second, we
show that this edge satisfies $2|A|+(|B|+|C|)/2>\lambda$, so that
integrality fixes the page sizes and excludes all other saturating
non-edges. Third, we prove that the Perron weight missing from $A$ is
so small that $A$ is completely joined to $U$. Finally, we complete
the proof by showing that no vertices remain outside $A\cup U$, so
that $G$ itself is a join.

Suppose that the theorem fails. Then there are $K_4$-free Nosal graphs
$G$ with $\lambda=\lambda(G)\to\infty$ and $f(G)\le\phi(\lambda)$. A
component attaining $\lambda$ is again such a graph, so we may assume
that every $G$ is connected. We use the notation of
Sections~\ref{subsec-2.3} and~\ref{sec-3}, and all $o(1)$ terms are
taken as $\lambda\to\infty$.

\emph{Step 1: an edge near the ratio $4:1:1$.} Case (iii) of
Section~\ref{subsec-2.4} holds with $\tau=1/\lambda$. By
Lemma~\ref{lem-3.1}, some edge of $H$ satisfies
$s_{uv}<(r/2)g_{uv}$. Name its endpoints as in Lemma~\ref{lem-3.2}, let
$A,B,C$ be the pages of $u^*uv$, and take $a,b,c,h,w$ from that lemma.
Thus $b=|B|/\lambda$, $c=|C|/\lambda$ and $a\le|A|/\lambda$. Put
$k=4r+h-2$. By Lemmas~\ref{lem-2.6} and~\ref{lem-2.8}, we get 
\[
 a^2+b^2+c^2\le\frac29+o(1),\qquad
 r^2+\frac54(b^2+c^2)\le\frac{11}{24}+o(1).
\]
Hence $r\le1-h^2/4$ for large $\lambda$, and Lemma~\ref{lem-3.3} yields
$w<2r$, $k>0$ and $\gamma(w)<0$. We claim that
\[
 a,r\to\frac49,\qquad b,c\to\frac19,\qquad w\to0.
\]
Let $(a_0,b_0,c_0,r_0,w_0)$ be a limit point of $(a,b,c,r,w)$, with
corresponding values $h_0,k_0,\gamma_0$. We show that this point must
be $(4/9,1/9,1/9,4/9,0)$. The two bounds above pass to the limit, and
also $0\le w_0\le2r_0$, $k_0\ge0$, $\gamma_0\le0$ and
$2b_0w_0=(r_0-a_0)(2b_0+r_0)$. Moreover, $k>0$ implies $r_0\ge1/3$. If
$b_0=0$, then $a_0=r_0$, so $k_0=4a_0+c_0-2\le\sqrt{34/9}-2<0$, a
contradiction. So $b_0>0$, and Lemma~\ref{lem-3.4} applies to the
limit: if $k_0>0$, then $\gamma_0>0$, which is false. Thus $k_0=0$.
Now $\gamma_0=w_0(2(2+h_0)r_0-w_0)\le0$ forces $w_0=0$, so $a_0=r_0$
and $4a_0+b_0+c_0=2$. Now equality in the Cauchy--Schwarz inequality
$(4a_0+b_0+c_0)^2\le18(a_0^2+b_0^2+c_0^2)$ proves the claim.

\emph{Step 2: the page sizes.} In what follows, our task is to turn
this limit into exact page sizes. Put $\sigma=2-h-4a$, so that
$k=4(r-a)-\sigma$. We first show that $\sigma<0$ and
\[
 r-a\le\Bigl(\frac12+o(1)\Bigr)(-\sigma).
\]
The limit $r/(2b)\to2$ implies $w=(3+o(1))(r-a)$. The inequality
$\gamma(w)<0$ reads $k>2w-w^2/((2+h)r)$. Since $(2+h)r\to80/81$ and
$w\to0$, we obtain $k>(2-o(1))w$ if $w>0$, and $k>0$ if $w=0$. Both
claims now follow from $-\sigma=k-4(r-a)$.

Put $p=|A|$ and $\ell=2p+(|B|+|C|)/2$. Since $y(A)\le p$, the
inequality $\sigma<0$ means that $2\lambda<4y(A)+|B|+|C|\le2\ell$, so
$\lambda<\ell$. By Lemma~\ref{lem-2.3}, the number
$f_0=\binom p2+\binom{|B|}2+\binom{|C|}2$ satisfies
$f_0\le f(G)\le\phi(\lambda)$. Therefore, Lemma~\ref{lem-4.1} applies:
\[
 |B|=|C|=d,\qquad p\in\{4d-2,4d-1\},\qquad f(G)=f_0,
\]
and $\varepsilon=\ell-\lambda>0$ satisfies
$\varepsilon(\ell+\lambda-3)\le1$. In particular, we have 
$\lambda\varepsilon\le1/2+o(1)$. Since $d/\lambda\to1/9$ by Step~1, we
also have $p/\lambda\to4/9$.

\emph{Step 3: complete adjacency.} From the above discussion, our
problem reduces to recovering a join from the eigenvalue equations.
Put $\eta=p-y(A)\ge0$ and $\omega_A=y(W\setminus A)=\lambda(r-a)\ge0$.
As $\ell=2p+d$, we have $\lambda\sigma=4\eta-2\varepsilon$, so Step~2
yields
\[
 \eta<\frac\varepsilon2,\qquad
 \omega_A\le(1+o(1))(\varepsilon-2\eta),\qquad
 \lambda\omega_A\le\frac12+o(1).
\]
In particular, $\eta=o(1)$, and every coordinate on $A$ is at least
$1-\eta$. If a vertex $v \in W\setminus A$ had a neighbor in $A$, its
eigenvalue equation would imply $\lambda\omega_A\ge\lambda y_v\ge1-\eta$,
a contradiction. Thus, no edge joins $A$ to $W\setminus A$, and every
vertex of $A$ has all its neighbors in $U$.

It remains to show that each vertex of $U$ is adjacent to all of $A$.
Subtracting the sum of the equations on $A$ from $p$ times the equation
at $u^*$, we get
\[
 \lambda\eta=\sum_{v\in U}(p-d_A(v))y_v.
\]
Every $v\in U$ has $d_A(v)\ge1$, because $u^*\in A$. If $d_A(v)<p$,
then the equation at $v$ implies $y_v\ge d_A(v)(1-\eta)/\lambda$, and
therefore
\[
 \lambda\eta\ge\frac{d_A(v)(p-d_A(v))(1-\eta)}\lambda
 \ge\frac{(p-1)(1-\eta)}\lambda\longrightarrow\frac49.
\]
But $\lambda\eta<\lambda\varepsilon/2\le1/4+o(1)$. Hence every vertex
of $U$ is adjacent to all of $A$, and $\eta=0$.

\emph{Step 4: no remaining vertices.} Combining the above facts, we
see that it remains to show that $W\setminus A$ is empty. By Step~2,
every saturating non-edge lies inside $A$, $B$ or $C$, so
Lemma~\ref{lem-4.2} implies $\Delta(H)\le d$ and
$e(H)\le d^2+|U|/2$. Also $p|U|\le m<\lambda^2$, so $|U|=O(\lambda)$.
Let $\beta$ be the largest coordinate at a vertex of positive degree
in $H$. The equation at a vertex attaining it reads
$\lambda\beta\le p+d\beta+\omega_A$. As $2p+d=\lambda+\varepsilon$,
this means that $2\beta-1\le(\varepsilon+2\omega_A)/(\lambda-d)$. Since
$\varepsilon+2\omega_A\le(3/2+o(1))/\lambda$, we obtain
\[
 \mu\le(2\beta-1)e(H)
 \le\frac{\varepsilon+2\omega_A}{\lambda-d}
          \left(d^2+\frac{|U|}{2}\right)
 \le\frac1{48}+o(1).
\]
If $W\setminus A$ were nonempty, connectedness would provide an edge
between it and $U$. Then $\nu\ge1-\omega_A=1-o(1)$, contradicting
$\mu>\nu$. Thus $W\setminus A$ is empty, and $G=I_p\vee H$. We can now
invoke Proposition~\ref{prop-4.3}: it states that
$f(G)>\phi(\lambda)$, which is the final contradiction.
\end{proof}

\section{The extremal family and the upper bound}\label{sec-5}

This section proves the upper bound in Theorem~\ref{thm-1.3} and the
sharpness of Theorem~\ref{thm-1.1}. We describe a family of $K_4$-free
graphs on the boundary $\lambda^2=m$ of the Nosal condition with
$\tF=\lambda^2/9$, perturb it into Nosal graphs, and use a matching to
reach each sufficiently large edge count.

\subsection{A family on the boundary
\texorpdfstring{$\lambda^2=m$}{lambda squared=m}}\label{subsec-5.1}

Let $t\ge1$ be an integer, and let $p\in(0,1/2]$ be a rational number with
$pt\in\mathbb Z$. Define the graph $G_{t,p}$ on four disjoint independent sets
$S_0,S_1,S_2,T$ with
\[
 |S_0|=(4+2p)t,\qquad |S_1|=|S_2|=pt,\qquad |T|=t .
\]
Join every vertex of $T$ to every vertex of $S=S_0\cup S_1\cup S_2$,
and join $S_1$ completely to $S_2$. These are all the edges of
$G_{t,p}$. The graph $G_{t,p}$ in the $p=1/4$ case was drawn in Figure~\ref{fig-1}(b). 

\begin{lemma}\label{lem-5.1}
The graph $G_{t,p}$ is $K_4$-free and satisfies
\begin{align*}
 e(G_{t,p})&=(2+p)^2t^2,\\
 f(G_{t,p})&=\Bigl(\frac12+p^2\Bigr)t^2-\Bigl(\frac12+p\Bigr)t,\\
 \lambda(G_{t,p})&=(2+p)t=\sqrt{e(G_{t,p})}.
\end{align*}
\end{lemma}
\begin{proof}
The graph is the join of $T$ and a triangle-free graph, namely
$K_{pt,pt}$ together with $(4+2p)t$ isolated vertices. Hence it is
$K_4$-free, and $e(G_{t,p})=(4+4p)t^2+(pt)^2=(2+p)^2t^2$. Exactly the pairs inside
$T$, $S_1$ and $S_2$ are saturating: their common neighborhoods
contain an edge, whereas every other non-edge has the independent
common neighborhood $T$. Thus
\[
 f(G_{t,p})=\binom t2+2\binom{pt}{2}
 =\left(\frac12+p^2\right)t^2-\left(\frac12+p\right)t.
\]
Finally, give the vertices of $S_0,S_1,S_2,T$ the values
$1/(2+p),1/2,1/2,1$, respectively. The resulting positive vector is
an eigenvector with eigenvalue $(2+p)t$, so this number equals the
spectral radius of the graph by the Perron--Frobenius theorem.
\end{proof}

The vertices of $T\cup S_1\cup S_2$, and only these vertices, lie in
triangles. Hence $q(G_{t,p})=(1+2p)t$, and it follows that
\[
 \frac{\tF(G_{t,p})}{\lambda(G_{t,p})^2}
 =\frac{1/2+p^2}{(2+p)^2}
 =\frac19+\frac{(4p-1)^2}{18(2+p)^2}.
\]
The minimum is therefore $1/9$, attained exactly when $p=1/4$. In
particular, for $4\mid t$, we get 
\[
 f(G_{t,1/4})=\frac9{16}t^2-\frac34t,\qquad
 q(G_{t,1/4})=\frac32t,\qquad
 \tF(G_{t,1/4})=\frac{\lambda(G_{t,1/4})^2}{9}.
\]

\subsection{Perturbation and exact edge counts}\label{subsec-5.2}

For the case $4\mid t$, we delete one vertex of $S_0$ from $G_{t,1/4}$ to obtain a
graph $G'$ with the same $f$ and $q$. Put $d=t/4$, $\ell=2t+d$ and
$\lambda=\lambda(G')$. Then $d<\lambda<\ell$, and the Perron equation
gives $\lambda^2-e(G')=d^2(\ell-\lambda)/(\lambda-d)>0$. Since
$e(G')=\ell^2-t$, it follows that $\tF(G')/\lambda^2\to1/9$. This proves that
the constant $1/9$ in Theorem~\ref{thm-1.1} is best possible. Also
$f(G')=e(G')/9-23t/36$ and $\sqrt{e(G')}<9t/4$, so the $\sqrt m$ error
order in Theorem~\ref{thm-1.3} is optimal.

\begin{proof}[Proof of the upper bound in Theorem~\ref{thm-1.3}]
We reach every edge count by a small change of $G_{t,1/4}$, namely a
join $I_a\vee(K_{d,d}\cup pK_2\cup sK_1)$ with $a$ close to $4d$ and
$2a+d$ slightly larger than $\sqrt m$. The idea behind the construction
is that deleting isolated vertices removes $a$ edges at a time, and a
matching provides the fine adjustment without changing $f$. The edge
count and the value of $f$ are then direct, and the point that needs
care is the strict Nosal condition.

\emph{The construction.} Let $m\ge500$ be an integer. Put $x=\sqrt m$, $\ell=\lceil x\rceil+4$
and $\varepsilon=\ell-x\in[4,5)$. Choose an integer $d\equiv\ell\pmod2$
nearest to $\ell/9$, and put $a=(\ell-d)/2$ and $j=a-4d$. Since
$\ell-9d$ is even and has absolute value at most $9$, the integer $j$
satisfies $|j|\le4$. Also $\ell\ge27$, and $d=2$ would require an even
$\ell\le27$, so $d\ge3$ and $a\ge4d-4\ge8$. Put
$g=\ell^2-m=\varepsilon(2\ell-\varepsilon)$, and note that
$\varepsilon\in[4,5)$ implies $8\ell-16\le g<10\ell-25$.

Start from $I_a\vee(K_{d,d}\cup(4a+2d)K_1)$, which has $\ell^2$ edges.
To remove exactly $g$ edges, set
\[
 k=\left\lceil\frac ga\right\rceil,\qquad
 p=ka-g,\qquad s=4a+2d-k-2p.
\]
Here $0\le p<a$. Since $a=(4\ell+j)/9$ and $|j|\le4$, it follows that
$g<10\ell-25<23a$, so $k\le23$. Consequently, we get 
$s\ge2a+2d-21\ge10d-29\ge1$, and the graph
\[
 G=I_a\vee\bigl(K_{d,d}\cup pK_2\cup sK_1\bigr)
\]
is defined and $K_4$-free. Its edge count is exactly
$a(2d+2p+s)+d^2+p=\ell^2-ka+p=m$. Only the pairs inside $I_a$ and the
two parts of $K_{d,d}$ are saturating, so by the identity in
Lemma~\ref{lem-4.1},
\[
 f(G)=\binom a2+2\binom d2=\phi(\ell)+\frac{j(j+3)}{18}.
\]

\emph{The Nosal condition.} The Perron vector is constant on the
orbits of the automorphism group, so the eigenvalue equations reduce to
\[
 \lambda=a\left(\frac{2d}{\lambda-d}
                 +\frac{2p}{\lambda-1}+\frac s\lambda\right),
 \qquad \lambda>d.
\]
Here $\lambda>d$, since $G$ is connected and properly contains
$K_{d,d}$. Also $x>d$, because $x-d=2a-\varepsilon>0$. The right side
decreases on $(d,\infty)$, so $\lambda(G)>x$ exactly when the right
side exceeds $x$ at $\lambda=x$. To rewrite this condition, a little
manipulation is required. Multiply by $x$, and use
$m=a(2d+2p+s)+d^2+p$. The condition becomes
$2ad^2/(x-d)+2ap/(x-1)>d^2+p$. Substituting $2a=x+\varepsilon-d$ turns
it into
\[
 d^2\varepsilon(x-1)>p(d-1-\varepsilon)(x-d).
\]
For $d\le5$ this is immediate, since $\varepsilon\ge4$. Let $d\ge6$.
Then $p\le a-1$, $x-d=2a-\varepsilon\le2a-4\le8d+4$ and
$d-1-\varepsilon\le d-5$. So the right side is less than
$8d^2(a-1)$. The left side is at least $4d^2(x-1)>8ad^2$, because
$x-1=2a+d-\varepsilon-1>2a$. Therefore $\lambda(G)>x=\sqrt m$.

\emph{The final count.} Since $j(j+3)/18\le14/9$ and
$\phi(\ell)-\ell^2/9=-\ell/3$, we conclude that
\[
 f(G)-\frac m9
 <\frac{10\ell-25}{9}-\frac\ell3+\frac{14}9
 =\frac{7\ell-11}{9}
 <\frac79x+\frac83<x,
\]
where the last inequality uses $x>12$. Consequently, we get 
$\psi(m)\le f(G)<m/9+\sqrt m$. Discarding isolated vertices bounds the
order by $2m$, so the minimum defining $\psi(m)$ is attained.
\end{proof}

\section{Concluding remarks}\label{sec-6}

This paper determines how many
$K_4$-saturating non-edges the Nosal
condition forces.
Theorem~\ref{thm-1.1} gives the sharp bound
$f(G)+\frac12q(G)>\frac19\lambda(G)^2$,
which extends Nosal's theorem to
$K_4$-free graphs.
Theorem~\ref{thm-1.3} shows that the jump
at the Nosal threshold has size
$\frac19m+O(\sqrt m\,)$ with explicit
constants, and Theorem~\ref{thm-1.4}
sharpens the error term at large spectral
radius. The constant $\frac19$ is less
than the constant $\frac8{33}$
that~\eqref{e-1} gives in terms of $m$ at
the Mantel threshold. Moreover, the extremal
graphs under the Nosal condition satisfy neither the Mantel condition
nor its spectral counterpart. Thus, the two
thresholds differ in both the sharp constant and
the extremal structure.

\smallskip 
The proof replaces a single maximizing edge
by a zero-sum comparison over all edges of a
neighborhood, and it compares the three
pages of a triangle on every edge with the
same coefficient. The Balogh--Liu theorem \cite{BL2014}
enters only through a weighted blow-up. 
We close with four related directions that these
methods suggest: the second-order term,
the range further above the Nosal threshold,
larger cliques and the corresponding problem at the spectral
Mantel threshold.

\medskip 
\noindent 
{\bf The second-order term in $\bm{\psi(m)}$.} 
The boundary graphs $G_{t,1/4}$ satisfy $f=\lambda^2/9-\lambda/3$
exactly, and Theorem~\ref{thm-1.4} gives the corresponding strict
inequality for Nosal graphs of large spectral radius. By
Theorem~\ref{thm-1.1}, the same inequality holds for every $K_4$-free
Nosal graph with $q(G)\le2\lambda(G)/3$, and $G_{t,1/4}$ satisfies
$q=2\lambda/3$. Since $x^2/9-x/3$ increases for $x\ge3/2$,
Theorem~\ref{thm-1.4} gives $\psi(m)>m/9-\sqrt m/3$ for sufficiently
large $m$. In the other direction, the Nosal graphs $G'$ of
Section~\ref{subsec-5.2} satisfy
$f(G')=m/9-23t/36<m/9-(23/81)\sqrt m$, where $m=e(G')$. Thus, the
lower limit of $(\psi(m)-m/9)/\sqrt m$ lies between $-1/3$ and
$-23/81$. It would be interesting to find its exact value.

\medskip 
\noindent 
{\bf Further above the Nosal threshold.}
In fact, Balogh and Liu \cite[Remark~(ii)]{BL2014} observed an interesting phenomenon that
$f(n, \lfloor n^2/4\rfloor + k) =\frac{2}{33}n^2 + O (n)$ for every
$1\le k\le \frac{n}{66}$. In other words, the minimum number of
$K_4$-saturating non-edges has a plateau: it stays close to
$\frac{2}{33}n^2$ even with linearly many additional edges above
$\lfloor n^2/4\rfloor$. We mention here that the triangle counts behave quite differently. By the
Lov\'{a}sz--Simonovits theorem~\cite{LS1975}, the minimum number of
triangles is $k\lfloor n/2\rfloor$ for every $1\le k< \frac{n}{2}$,
so it grows by $\lfloor n/2\rfloor$ with each additional edge. Wang,
Yang and Zheng~\cite{WYZ2026} later determined the asymptotic minimum
of $f(n,m)$ for every $m$ up to $\mathrm{ex}(n,K_4)$.

\smallskip 
In the spectral setting, Chen, Li and Tang~\cite{CLT2026}
proved that $\lambda (G) \ge \sqrt{m} +c$ forces more than $c m$
triangles for every real $c >0$. This can be viewed as a spectral
counterpart of the Lov\'{a}sz--Simonovits theorem. The analogue for
$K_4$-saturating non-edges is the function
\[
 \psi(m, c) :=\min\bigl\{f(G): e(G)=m,\ G\text{ is }
 K_4\text{-free},\ \lambda(G)\ge \sqrt{m} +c \bigr\}.
\]
Fix $c>0$, and let $m$ be sufficiently large. Since $\phi$
increases on $[3/2,\infty)$, Theorem~\ref{thm-1.4} gives the lower
bound below. For the upper bound, take the graph
$G=I_a\vee(K_{d,d}\cup pK_2\cup sK_1)$ from the proof of the upper
bound in Theorem~\ref{thm-1.3}, but with
$\ell=\lceil\sqrt m+145c\rceil+4$.
Put $x=\sqrt m$ and $\varepsilon=\ell-x$, so that
$145c+4\le\varepsilon<145c+5$. Then $k=\lceil g/a\rceil$ is
$O(c+1)$, while $a$ and $d$ grow linearly in $x$. Hence
$s\ge1$ for large $m$, and the graph is well defined.
 The computation in that
proof shows that $\lambda(G)\ge x+c$ if and only if
\[
 \frac{d^2(\varepsilon-c)}{x+c-d}
 +\frac{p(\varepsilon-c-d+1)}{x+c-1}\ge 2cx+c^2.
\]
Since $d=(\frac19+o(1))x$, $p<a=(\frac49+o(1))x$ and
$\varepsilon\ge145c+4$, the left side is at least
$(2c+\frac1{162}+o(1))x$, so this holds. The identity
$f(G)-\frac m9=\frac g9-\frac\ell3+\frac{j(j+3)}{18}$ from
that proof, with $g=2\varepsilon x+\varepsilon^2$, then
gives 
\[
 \frac19m+\frac{2c-3}{9}\sqrt m+\frac19 c^2-\frac13 c
 <\psi(m,c)
 \le\frac19m+\frac{290c+7}{9}\sqrt m+O(c^2 +1).
\]
Hence, the first-order term $\frac19m$ does not depend on $c$,
whereas the triangle bound $cm$ grows linearly with $c$. So the
plateau also appears at the Nosal threshold. The first-order term
does change when $c$ has order $\sqrt m$, since
$\lambda(G)\ge\beta\sqrt m$ with $\beta>1$ forces
$f(G)>\frac19\beta^2m-O(\sqrt m\,)$ by Theorem~\ref{thm-1.3}. 
Two natural questions remain. What is the true coefficient of
$\sqrt m$ in $\psi(m,c)$? For fixed
$\beta\in(1,2/\sqrt3]$, what is the asymptotic minimum of
$f(G)$ over all $K_4$-free graphs with $m$ edges and
$\lambda(G)\ge\beta\sqrt m$? The bound
$\frac19\beta^2m$ is far from sharp for large $\beta$.
For example, the regular Tur\'{a}n graph $T_{n,3}$
satisfies $\beta=2/\sqrt3$ and $f (T_{n,3})= \big(\frac12+o(1)
\big)m$, which is much larger than $\frac{4}{27}m$ of
this bound. 

\medskip 
\noindent
{\bf Saturating non-edges for larger cliques.}  
Nikiforov's inequality \cite{Nik2002} 
shows that every graph $G$ with $m$ edges and 
$\lambda(G)>\sqrt{2m(1-1/r)}$ contains a copy of
$K_{r+1}$. The corresponding saturation problem
concerns $K_{r+2}$-free graphs $G$ with $m$ edges and
$\lambda(G)>\sqrt{2m(1-1/r)}$. It asks for the minimum of
$f_{r+2}(G)$, the number of $K_{r+2}$-saturating non-edges of $G$. For
$r=2$, Theorem~\ref{thm-1.3} shows that this minimum is
$(1/9+o(1))m$. 
The transference of Section~\ref{subsec-2.2} extends to every $r$.
He, Ma, Ma and Ye~\cite{HMMY2023} determined the minimum number of
$K_{r+2}$-saturating non-edges in $K_{r+2}$-free graphs with $n$
vertices and $\mathrm{ex}(n,K_{r+1})+1$ edges. This minimum equals
$(c_r+o(1))n^2$, where
\[
 c_r=\frac{2(r-1)^2}{(r+1)(4r^2-3r+1)},\qquad c_2=\frac2{33}.
\]
Let $G$ be a connected $K_{r+2}$-free graph with
$\lambda(G)> \sqrt{2m(1-1/r)}$ and Perron vector $\by$. The proof of
Lemma~\ref{lem-2.2}, with their theorem \cite{HMMY2023} in place of~\eqref{e-1}, gives
\[
 \sum_{uv}y_u^2y_v^2+
 \frac12\sum_{v}y_v^4
 \  \ge \ c_r\|\by\|_2^4,
\]
where the first sum runs over the $K_{r+2}$-saturating non-edges of
$G$ and the second sum runs over the vertices in copies of $K_{r+1}$. Indeed,
the weights $z_v=y_v^2/\|\by\|_2^2$ satisfy
$\sum_{uv\in E(G)}z_uz_v\ge\lambda(G)^2/(4m)>(1-1/r)/2$, so the blow-up
has more than $\mathrm{ex}(s,K_{r+1})+1$ edges for large $s$.
Moreover, a pair inside the class of $v$ is saturating exactly when
$v$ lies in a copy of $K_{r+1}$. This inequality provides the first
stage of our argument for every $r\ge3$. The second stage, an analogue
of the edge comparison of Section~\ref{sec-3}, remains open, and  
some new ideas are needed. 

\smallskip 
Complete $r$-partite graphs with equal parts lie on the
boundary $\lambda(G)^2=2m(1-1/r)$. A natural analogue of
$G_{t,p}$ embeds a complete bipartite graph into one part
of a complete $r$-partite graph. A numerical optimization
of its part sizes on the boundary gives
$f\approx0.1522\lambda^2$ for $r=3$, 
$f\approx0.1201\lambda^2$ for $r=4$ 
and  $f\approx 0.0977\lambda^2$ for $r=5$. 
On the other hand, 
the edge-extremal
graphs of He, Ma, Ma and Ye~\cite{HMMY2023} satisfy
$\lambda\ge 2m/n$, so they have
$f\le(c_r+o(1))(\frac r{r-1})^2\lambda^2$. This
coefficient equals $0.1607$, $0.1208$ and $0.0969$ for
$r=3,4,5$, respectively. 
So the analogue of $G_{t,p}$ may be better
only for $r=3 ,4$, but the edge-extremal graphs are better for every 
$r\ge 5$. It is possible that the extremal graphs for larger
cliques compete and combine features of both families.

\medskip 
\noindent 
{\bf Above the spectral Mantel threshold.} 
The condition $\lambda(G)>\lambda(T_{n,2})$
 is the spectral form of the Mantel condition $e(G)> e(T_{n,2})$
for graphs of order $n$. 
At this threshold, many spectral extremal problems have the
same extremal graphs as their edge versions. Examples include
cliques \cite{Wil1986,Nik2007}, odd cycles \cite{Nik2008-odd-cycles}, books \cite{ZL2022jgt} and
color-critical graphs \cite{Nik2009-critical}.
In general, Wang, Kang and Xue \cite{WKX2023} proved that
for sufficiently large $n$, every spectral extremal graph for $F$ is
edge-extremal whenever $\mathrm{ex}(n,F)=e(T_{n,r})+O(1)$. We
refer the reader to~\cite{FTZ2025,Byr2024,FGKLZ2027} for related
extensions.

\smallskip  
We next illustrate that under the condition $\lambda(G)>\lambda(T_{n,2})$, 
 the minimum of $f(G)$ over all $K_4$-free graphs on $n$ vertices is
 smaller than that in~\eqref{e-1}. The family $G_{t,p}$ lies
below this spectral threshold, since
$\lambda(G_{t,p})/|V(G_{t,p})|=(2+p)/(5+4p)<1/2$. Recall from
Figure~\ref{fig-1}(a) that the Balogh--Liu graph is a blow-up of the
$5$-cycle $V_1V_2V_3V_4V_5$ with the chord $V_1V_3$. It satisfies the
spectral Mantel condition since its size exceeds $\lfloor n^2/4\rfloor$. 
However, its class sizes are not optimal under the spectral Mantel condition. We reweight the classes of the same blow-up and take
\[
{|V_1|} = {|V_3|} =\Big(\tfrac{1}{48} \sqrt6 \Big)^{\!1/2} n ,
 \quad
{|V_2|}=\Bigl(\tfrac{1}{24} \big(5\sqrt6-12 \big)\Bigr)^{\!1/2}n,
 \quad |V_4|=|V_5|=\tfrac12\bigl(n-|V_1|-|V_2|-|V_3|\bigr).
\]
By computation, the resulting graphs satisfy $\lambda(G)=n/2+O(1)$. To
make the inequality strict, note that $N(V_4)\subsetneq N(V_1)$ in the
underlying five-vertex graph. The eigenvalue equations therefore show
that the Perron vector of $G$ takes a larger value $z_1$ on $V_1$ than
its value $z_4$ on $V_4$. Move $\varepsilon n$ vertices from $V_4$ to
$V_1$, and give them the value $z_1$. Since $V_1$ and $V_4$ are
nonadjacent in the underlying graph, this raises the quadratic form by
$2\lambda(G)(z_1^2-z_4^2)\varepsilon n$ and the squared norm by
$(z_1^2-z_4^2)\varepsilon n$. Thus, the Rayleigh principle shows that 
the spectral radius increases by a
positive multiple of $\varepsilon n$. The new graph is still a blow-up of the same
five-vertex graph, and its saturating non-edges are the
pairs inside $V_1$, $V_2$ and $V_3$. So $f(G)$ changes by
$O(\varepsilon n^2)$. Letting $\varepsilon\to0$ with
$\varepsilon n\to\infty$, we obtain graphs that satisfy the spectral condition $\lambda (G) > n/2 \ge  \lambda (T_{n,2})$
and
\[
 f(G)= \Bigl(\tfrac{\sqrt6-2}{8}+o(1)\Bigr)n^2,\qquad
 \tfrac{\sqrt6-2}{8}=0.0561\cdots<\tfrac2{33}.
\]
Consequently, the spectral Mantel threshold $\lambda(G) > \lambda (T_{n,2})$ reveals a different 
phenomenon and 
 forces fewer $K_4$-saturating non-edges than the classical Mantel threshold $e(G)> e(T_{n,2})$. 

\smallskip  
Theorem~\ref{thm-1.3} also gives an explicit lower bound at the spectral
Mantel threshold. Indeed, let
$G$ be a $K_4$-free graph on $n$ vertices with
$\lambda(G)>\lambda(T_{n,2})$. If $G$ is a Nosal graph,
then Theorem~\ref{thm-1.3} gives
$f(G)>\lambda(G)^2/9-\lambda(G)\ge n^2/36-O(n)$. Otherwise,
we have $e(G)\ge\lambda(G)^2>\lfloor n^2/4\rfloor$. Deleting edges
as in the proof of Lemma~\ref{lem-2.2} then gives
$f(G)\ge f(n,\lfloor n^2/4\rfloor+1)
= \big(\frac{2}{33}+o(1) \big)n^2$. Hence, the minimum of $f(G)$ at
the spectral Mantel threshold lies between $\big(\frac1{36}-o(1) \big)n^2$ and
$\big(\frac{1}{8}(\sqrt6-2)+o(1) \big)n^2$. We conjecture that the upper bound is
sharp. Determining the correct constant together with the
corresponding extremal graphs remains open. For larger cliques, we
expect that the extremal graphs at the spectral threshold
$\lambda(G)>\lambda(T_{n,r})$ are blow-ups of the same underlying
graphs as in He, Ma, Ma and Ye~\cite{HMMY2023}, although the part sizes may differ.

\section*{Acknowledgments}
The author would like to thank Hong Liu for 
introducing him to the spectral problem of counting $K_4$-saturating non-edges 
during the first IBS
ECOPRO student research program in fall 2023. 
This work was supported by the China Postdoctoral Science Foundation (No.~2026M793384), the Postdoctoral Fellowship Program of the China Postdoctoral Science Foundation (No.~GZC20261656) and the National Natural Science Foundation of China (No.~12601675). 
During the preparation of this paper, the author used the AI tools to brainstorm proof strategies and to provide simplifications of several proof arguments. 
All mathematical arguments were written, reviewed and verified by the author, who takes full responsibility for the content of this paper.

\end{document}